\documentclass[review]{elsarticle}
\usepackage{amsthm,amsmath,amssymb,amsopn}
\usepackage{graphicx,epstopdf,epsfig,subfig}
\usepackage[colorlinks,linktocpage,linkcolor=blue]{hyperref}
\usepackage[table,xcdraw]{xcolor}
\usepackage{url,color}
\usepackage{booktabs}
\usepackage{colortbl,multirow}
\usepackage{tikz}
\usetikzlibrary{shapes.geometric}
\usepackage[T1]{fontenc}
\usepackage{algorithm,algorithmic}
\numberwithin{equation}{section}
\numberwithin{table}{section}

\theoremstyle{definition}

\def\dr{\mathrm{d}r}
\def\dz{\mathrm{d}z}

\def\det{\text{det}}

\newcommand{\gv}[1]{\ensuremath{\boldsymbol{#1}}}
\newcommand{\grad}[1]{\gv{\nabla} #1} 
\renewcommand{\div}[1]{\gv{\nabla} \cdot #1} 

\definecolor{lightgray}{gray}{0.9}

\journal{Journal of Computational Physics}

\begin{document}
\begin{frontmatter}
\title{High-order  staggered Lagrangian hydrodynamics (II) : the artificial viscosity and hourglass control algorithm}
\author[add1]{Zhiyuan Sun}
\ead{zysun.math@gmail.com}

\author[add1]{Jun Liu}
\ead{caepcfd@126.com}

\author[add1]{Pei Wang}
\ead{wangpei@iapcm.ac.cn}

\address[add1]{Institute of Applied Physics and Computational Mathematics, Beijing 100094, P. R. China}

\begin{abstract}
In this article, we investigate the artificial viscosity and hourglass control algorithms for high-order staggered Lagrangian hydrodynamics(SGH), as proposed in~\cite[Sun et al., 2025]{Sun2025High}. Inspired by the subzonal pressure method in classical staggered Lagrangian hydrodynamics, we extend the stiffness-based hourglass control algorithm to the high-order setting, enriching the pressure field from the $Q^{m-1}$ to the $Q^{m}$ polynomial space. A unified framework for this hourglass control approach is established, from which the classical subzonal pressure method naturally emerges as the special case of the $Q^1-P^0$ space. The artificial viscosity follows the formulation in~\cite[Dobrev et al., 2012]{Dobrev2012High}. We show that the viscosity admits a concise form, with intermediate variables explicitly computable, leading to improved computational efficiency and easier implementation. Moreover, the tensor viscosity in classical staggered Lagrangian hydrodynamics can be derived in a similarly compact and explicit form. Numerical experiments on two-dimensional problems are presented to demonstrate the accuracy and efficiency of the proposed algorithms.

\noindent\textbf{keyword:} high-order, Staggered Lagrangian hydrodynamics, article viscosity, hourglass control
  
  
\end{abstract}


\end{frontmatter}
\section{Introduction}\label{sec:introduction}
This paper is a continuation of our previous work on high-order staggered Lagrangian hydrodynamics (SGH) presented in~\cite{Sun2025High}. We focus on two key algorithms in high-order SGH: the hourglass control algorithm and the artificial viscosity algorithm. In classical SGH, mesh tangling poses a principal challenge in two and three dimensions~\cite{margolin2013arbitrary}. Beyond the physical motion of the fluid, hourglass distortion\textemdash an artifact of mesh tangling caused by deficiencies in the numerical scheme~\cite{Menchen1964tensor}\textemdash is a particularly critical issue. The hourglass control algorithm is essential for classical SGH and generally follows two strategies: the viscosity-based method~\cite{Flanagan1981uniform,Margolin1987method} and the stiffness-based method~\cite{Caramana1998Elimination,Cheng2015Elimination,Sun2022Understanding}. However, the high-order curvilinear finite element method for Lagrangian hydrodynamics proposed in~\cite{Dobrev2011Curvilinear,Dobrev2012High} claims that no hourglass filter is required due to the strong mass conservation principle~\cite{Anderson2020Multiphysics}. In our previous work~\cite{Sun2025High}, we demonstrated that enforcing strong mass conservation with sufficiently accurate quadrature rules leads to inconsistencies in the thermodynamic variables, while ensuring thermodynamic consistency inevitably introduces hourglass distortion.

This paper extends the methodology of our previous work~\cite{Sun2022Understanding}, adapting the classical SGH method for high-order applications. For the $Q^m-Q^{m-1}$ finite element space pair in high-order SGH, the kinematic variables are approximated in the continuous $Q^m$ space, while the thermodynamic variables are approximated in the discontinuous $Q^{m-1}$ space. In our construction of the hourglass control algorithm, the pressure term is enriched from the $Q^{m-1}$ polynomial space to the $Q^m$ polynomial space. According to numerical quadrature theory, the $Q^{m-1}$ pressure term can be integrated exactly using an $m^2$-point Gauss-Legendre quadrature rule, which is equivalent to a $(m+1)^2$-point rule for the same integrand. The pressure variation evaluated with the $(m+1)^2$-point rule can be interpreted as the anti-hourglass force. Consistent with the $Q^1-P^0$ case, the pressure variation in the high-order setting is predicted by the product of density variation and the square of the sound speed. The high-order density term arises naturally from the $(m+1)^2$-point quadrature rule, analogous to the subzonal mass conservation in~\cite{Caramana1998Elimination}. With the defined interpolation matrix, density variation can be computed point wise. The classical subzonal pressure method~\cite{Caramana1998Elimination} is recovered as a special case of the proposed framework when $m=1$.

Artificial viscosity is a widely used numerical technique for capturing shock waves and suppressing oscillations in SGH simulations. Extensive work has been devoted to this topic; we refer the reader to the review articles~\cite{Benson1992Computational,Barlow2016Arbitrary} for an overview. In this study, we adopt the form of artificial viscosity proposed in~\cite{Kolev2009tensor,Dobrev2012High}. Rather than advancing the theory or introducing a novel formulation, our focus is on improving the computational efficiency of this viscosity. In classical SGH, the implementation of bulk viscosity~\cite{Schulz1964Two,Wilkins1980Use} and tensor viscosity~\cite{Campbell2001tensor} is both efficient and economical in terms of computational resources. The canonical open-source hydrocode DYNA2D~\cite{Hallquist1980User} is particularly well known for its efficiency. In such hydrocodes, the original mathematical formulas are decomposed into intermediate variables that must be computed and reused wherever possible. However, in the high-order case, fully decomposing the discretized formulas requires a large number of intermediate variables, which makes programming cumbersome.

To balance efficiency with ease of implementation, we introduce a concise formulation of tensor viscosity. This formulation allows the classical tensor viscosity~\cite{Campbell2001tensor} to be derived conveniently. Using the $Q^1-P^0$ case as an example, we present a streamlined code implementation of classical tensor viscosity. Although this approach incurs a slightly higher computational cost compared to traditional formulas, it significantly simplifies the implementation, which is particularly valuable in the high-order setting. Furthermore, we demonstrate that the viscosity parameters can be expressed explicitly, enabling an efficient practical implementation of artificial viscosity.

The roles of artificial viscosity and hourglass control are briefly introduced. Although they are distinct, the two are closely related algorithms. In this work, we address them simultaneously for three reasons. First, under the proposed framework, both algorithms rely on the same Gauss-Legendre quadrature rule with a given level of precision. As a result, each requires the gradients of the shape functions at the same quadrature points. Moreover, both algorithms demand interpolation of thermodynamic variables from $m^2$ points to $(m+1)^2$ points, using the same interpolation matrix. Second, artificial viscosity is also capable of suppressing hourglass modes. The main conclusion of our previous work~\cite{Sun2022On} for classical SGH extends naturally to the high-order case: high-order quadrature rules for artificial viscosity are beneficial for hourglass control, though they inevitably compromise robustness. Third, both algorithms directly affect the robustness of high-order SGH. The singularity of the Jacobian determinant is a primary source of instability. As the Jacobian determinant at a quadrature point approaches zero, the corresponding density variation tends to infinity, producing an inappropriate anti-hourglass force. Similarly, since the formulation of artificial viscosity contains the reciprocal of the Jacobian determinant, it also yields unreasonably large viscosity forces under such conditions. Therefore, for two-dimensional $Q^{m}-Q^{m-1}$ finite element space pairs, we specify the $(m+1)^2$-point quadrature rule for the discretization of both artificial viscosity and hourglass control. This choice provides a balance between robustness and computational efficiency.

This paper is arranged as follows. Section~\ref{sec:hourglass} details the construction of the anti-hourglass force, which is derived from the difference between two pressure terms. We assume a higher-order polynomial pressure term can be obtained and then introduce an interpolation matrix to compute the density variation point-wise. In Section~\ref{sec:artificial_viscosity}, we present a concise formulation of artificial viscosity and demonstrate its programming efficiency. We also show that the viscosity parameter can be explicitly calculated. The energy evolution with both the hourglass and artificial viscosity forces is discussed in Section~\ref{sec:energy_conservation}. Finally, Section~\ref{sec:numerical_results} presents numerical results for two-dimensional benchmark problems, illustrating the efficiency of the proposed algorithms.


\section{The hourglass control algorithm}\label{sec:hourglass}
We only consider the momentum equation of the compressible Euler equation in Lagrangian description.
\begin{equation}\label{eq:momentum}
\rho \frac{d \vec{u} }{d t} =- \nabla p  
\end{equation} 
Recall the notations in~\cite{Sun2025High}. The kinematic variables are discretized in the continuous $Q^m$ finite element space, while the thermodynamic variables are discretized in the discontinuous $Q^{m-1}$ finite element space. For an element $K$, let $N_i, N_j$ denote the shape functions of the $Q^m$ space for kinematic variables, and let $\phi_k, \phi_l$ denote the shape functions of the $Q^{m-1}$ space for thermodynamic variables. In each element, the number of degrees of freedom (DOFs) in the $Q^m$ space is denoted by $kdof$, and the number of DOFs in the $Q^{m-1}$ space is denoted by $tdof$.\begin{displaymath}
  \begin{aligned}
      \vec{u}_h|_K &=\sum_{i=1}^{kdof} \vec{u}_{i}N_{i},\\
      \rho_h|_K &=\sum_{k=1}^{tdof} \rho_{k}\phi_{k},\\
      p_h|_K &=\sum_{k=1}^{tdof} p_{k}\phi_{k},\\
      e_h|_K &=\sum_{k=1}^{tdof} e_{k}\phi_{k}.
\end{aligned}
\end{displaymath}

We now semi-discretize the momentum equation~\eqref{eq:momentum} using the Petrov-Galerkin technique.
\begin{equation}\label{eq:momentum_conservation}
\sum_{i=1}^{kdof} \frac{d  \vec{u}_{i}}{d t} \int_{K} \rho N_{i} N_j \dr\dz = \int_{K} p \nabla N_j \dr\dz,
\end{equation}
Only the velocity variable is expressed in discrete form. The mass matrix on the left-hand side (LHS) is evaluated using the Gauss-Lobatto quadrature rule, which yields a diagonal mass matrix. Since the mass matrix remains invariant during time evolution, its diagonal entries can be accumulated to the nodes to obtain the nodal mass. Once the right-hand side (RHS) is computed, the nodal forces are accumulated in the same manner; dividing them by the nodal mass then yields the nodal accelerations.

The actual pressure term $p_h$ is a $Q^{m-1}$ polynomial, whose number of DOFs is insufficient to constrain the hourglass modes. We now focus on the RHS of equation~\eqref{eq:momentum_conservation}. Since $p_h \in Q^{m-1}$, the RHS of equation~\eqref{eq:momentum_conservation} can be evaluated exactly using the $m^2$-point Gauss–Legendre quadrature rule:
\begin{equation}\label{p_h_RHS}
\int_{K} p_h \nabla N_j \dr\dz=  \sum_{q=1}^{m^2} \omega_{q} p_h(\gv{\xi}_q) \nabla N_j(\gv{\xi}_q) \det J(\gv{\xi}_q),
\end{equation}
where the highest polynomial degree of the integrand is $2m-1$, which can be integrated exactly with the $m^2$-point Gauss-Legendre quadrature. Here, $\gv{\xi}_q$ and $\omega_q$ denote the quadrature points and weights, respectively, and $J$ is the Jacobian matrix of the mapping from the reference element to the physical element.

If we assume that a higher-order pressure field $\tilde{p}$ can be obtained\textemdash that is, a pressure polynomial $\tilde{p} \in Q^{m}$  with additional DOFs\textemdash then replacing $p_h$ in equation~\eqref{eq:momentum_conservation} with $\tilde{p}$ yields a new RHS that eliminates the hourglass deficiencies. Now the highest polynomial degree of the integrand is $2m$, which can be integrated exactly using the $(m+1)^2$-point Gauss-Legendre quadrature rule:
\begin{equation}\label{p_tilde_RHS}
\int_{K} \tilde{p}\nabla N_j \dr\dz  = \sum_{q=1}^{(m+1)^2} \omega_{q}\tilde{p}(\gv{\xi}_q) \nabla N_j(\gv{\xi}_q) \det J(\gv{\xi}_q).
\end{equation}
It is important to note that the quadrature points $\gv{\xi}_q$ and weights $\omega_q$ are distinct for different quadrature rules and therefore cannot be directly combined through addition or subtraction. For notational simplicity, however, we employ the same symbols throughout this work.

The RHS of equation~\eqref{p_h_RHS} leads to hourglass distortion, whereas the RHS of equation~\eqref{p_tilde_RHS} suppresses it. A natural approach, therefore, is to construct the anti-hourglass force from the difference between these two RHS terms.
\begin{equation}\label{eq:hg_force}
\begin{split}
\vec{fh}_{j}&= \int_{K} \tilde{p} \nabla N_j \dr\dz - \int_{K} p_h \nabla N_j \dr\dz \\
      &=\sum_{q=1}^{(m+1)^2} \omega_{q} \tilde{p}(\gv{\xi}_q) \nabla N_j(\gv{\xi}_q) \det J(\gv{\xi}_q) - \sum_{q=1}^{m^2} \omega_{q}p_h(\gv{\xi}_q) \nabla N_j(\gv{\xi}_q) \det J(\gv{\xi}_q) \\
      &= \sum_{q=1}^{(m+1)^2} \omega_{q}\tilde{p}(\gv{\xi}_q) \nabla N_j(\gv{\xi}_q) \det J(\gv{\xi}_q) - \sum_{q=1}^{(m+1)^2} \omega_{q}p_h(\gv{\xi}_q) \nabla N_j(\gv{\xi}_q) \det J(\gv{\xi}_q)\\
      &=\sum_{q=1}^{(m+1)^2} \omega_{q} (\tilde{p}-p_h)(\gv{\xi}_q) \nabla N_j(\gv{\xi}_q) \det J(\gv{\xi}_q) \\
      &=\sum_{q=1}^{(m+1)^2} \omega_{q} \delta p(\gv{\xi}_q) \nabla N_j(\gv{\xi}_q) \det J(\gv{\xi}_q)
\end{split}
\end{equation}
where $\delta p$ denotes $\tilde{p}-p_h$. The third equality follows from the exactness of both the $m^2$-point and $(m+1)^2$-point Gauss-Legendre quadrature rules for the RHS of equation~\eqref{p_h_RHS}. The anti-hourglass force $\vec{fh}_j$ is evaluated using the $(m+1)^2$-point quadrature rule, and $\delta p$ must be determined at these points.

A natural choice for $\delta p$ is to approximate it as the product of the density variation and the square of the sound speed.
\begin{equation}\label{pressure_variation}
  \delta p= c_s^2 \delta \rho. 
\end{equation}
Here, $c_s$ denotes the sound speed, and $\delta \rho$ represents the density variation. Although these relationships are straightforward in the continuous setting, their discrete evaluation is more involved. In our framework, thermodynamic variables are defined at $m^2$ Gauss-Legendre points, whereas $\delta p$ is required at $(m+1)^2$ points. Focusing first on the sound speed $c_s$ and setting aside the density variation $\delta \rho$, the discrete sound speed $[c_s]_h$ is defined at the $m^2$ Gauss quadrature points. Accordingly, an interpolation from the $m^2$-point variable $[c_s]_h$ to the $(m+1)^2$-point variable $\tilde{c_s}$ is required, which can be conveniently expressed in vector form:
\begin{equation}
  \begin{bmatrix}
    \tilde{c_s}_1\\
    \tilde{c_s}_2\\
    ...\\
    \tilde{c_s}_{(m+1)^2}
\end{bmatrix}
=
M_{interp}  \begin{bmatrix}
  [c_s]_{h1}\\
  [c_s]_{h2}\\
  ...\\
  [c_s]_{hm^2}
\end{bmatrix},
\end{equation}
where $M_{\text{interp}}$ denotes the interpolation matrix, whose size is $(m+1)^2 \times m^2$ and depends on the FEM space pair. The vector $\tilde{c_s}$ represents the interpolated sound speed at the $(m+1)^2$ Gauss-Legendre quadrature points.

For the density variation $\delta \rho$, we draw on the concept of subzonal density from compatible hydrodynamics algorithms~\cite{Caramana1998Elimination}. Assuming that the mass conservation law is satisfied at the $(m+1)^2$ quadrature points, the high-order density variable $\tilde{\rho}$ can thus be obtained.
\begin{equation}\label{discrete_mass_conservation_mp1}
\begin{split}
\tilde{\rho}^{n}(\gv{\xi}_{1}) \det J^{n}(\gv{\xi}_{1})&=\rho_{0}(\gv{\xi}_{1}) \det J_{0}(\gv{\xi}_{1}),\\
\tilde{\rho}^{n}(\gv{\xi}_{2}) \det J^{n}(\gv{\xi}_{2})&=\rho_{0}(\gv{\xi}_{2}) \det J_{0}(\gv{\xi}_{2}),\\
\tilde{\rho}^{n}(\gv{\xi}_{i}) \det J^{n}(\gv{\xi}_{i})&=\rho_{0}(\gv{\xi}_{i}) \det J_{0}(\gv{\xi}_{i}), ...\\
\tilde{\rho}^{n}(\gv{\xi}_{(m+1)^2}) \det J^{n}(\gv{\xi}_{(m+1)^2})&=\rho_{0}(\gv{\xi}_{(m+1)^2}) \det J_{0}(\gv{\xi}_{(m+1)^2})
\end{split}
\end{equation}
The density variable $\rho_h$ is also defined at the $m^2$ Gauss-Legendre quadrature points and must be interpolated to the $(m+1)^2$ quadrature points. The density variation can then be calculated as follows:
\begin{equation}\label{desity_variation}
  \begin{bmatrix}
    \delta \rho_1\\
    \delta \rho_2\\
    ...\\
    \delta \rho_{(m+1)^2}
\end{bmatrix}  
=
\begin{bmatrix}
  \tilde{\rho}_1\\
  \tilde{\rho}_2\\
  ...\\
  \tilde{\rho}_{(m+1)^2}
\end{bmatrix}  
- M_{interp}
\begin{bmatrix}
  \rho_{h1}\\
  \rho_{h2}\\
  ...\\
  \rho_{hm^2}
\end{bmatrix}.
\end{equation}

The anti-hourglass force $\vec{fh}_j$ can now be computed using equation~\eqref{eq:hg_force}. We next discuss the implementation of this algorithm on a case-by-case basis.
\subsection*{$Q^1-P^0$ case}
The $Q^1-P^0$ pair is the traditional finite element space used in Lagrangian hydrodynamics. The main results are similar to those of the subzonal hourglass control algorithm~\cite{Sun2022Understanding}, with the only difference being the definition of the subzonal volume. In this case, the interpolation matrix reduces to a simple vector.
\begin{displaymath}
  M_{interp}=
  \begin{bmatrix}
    1\\
    1\\
    1\\
    1
  \end{bmatrix}.  
\end{displaymath}

The density variable $\tilde{\rho}$ is calculated from the mass conservation law and belongs to the $Q^1$ finite element space.
\begin{equation}\label{discrete_mass_conservation_Q1}
  \begin{split}
  \tilde{\rho}^{n}(-\frac{1}{\sqrt{3}},-\frac{1}{\sqrt{3}}) \det J^{n}(-\frac{1}{\sqrt{3}},-\frac{1}{\sqrt{3}})&=\rho_{0}(-\frac{1}{\sqrt{3}},-\frac{1}{\sqrt{3}}) \det J_{0}(-\frac{1}{\sqrt{3}},-\frac{1}{\sqrt{3}}),\\
  \tilde{\rho}^{n}(\frac{1}{\sqrt{3}},-\frac{1}{\sqrt{3}}) \det J^{n}(\frac{1}{\sqrt{3}},-\frac{1}{\sqrt{3}})&=\rho_{0}(\frac{1}{\sqrt{3}},-\frac{1}{\sqrt{3}}) \det J_{0}(\frac{1}{\sqrt{3}},-\frac{1}{\sqrt{3}}),\\
  \tilde{\rho}^{n}(-\frac{1}{\sqrt{3}},\frac{1}{\sqrt{3}}) \det J^{n}(-\frac{1}{\sqrt{3}},\frac{1}{\sqrt{3}})&=\rho_{0}(-\frac{1}{\sqrt{3}},\frac{1}{\sqrt{3}}) \det J_{0}(-\frac{1}{\sqrt{3}},\frac{1}{\sqrt{3}}),\\
  \tilde{\rho}^{n}(\frac{1}{\sqrt{3}},\frac{1}{\sqrt{3}}) \det J^{n}(\frac{1}{\sqrt{3}},\frac{1}{\sqrt{3}})&=\rho_{0}(\frac{1}{\sqrt{3}},\frac{1}{\sqrt{3}}) \det J_{0}(\frac{1}{\sqrt{3}},\frac{1}{\sqrt{3}}).
  \end{split}
\end{equation}
The only difference between the proposed algorithm and the subzonal hourglass control algorithm~\cite{Sun2022Understanding} lies in the volume calculation. In the latter, the subzonal volume is defined as the area enclosed by the barycenter of the grid, one node, and the two midpoints of the adjacent edges. In our approach, the Jacobian determinant is used directly as the subzonal volume. The computational procedure for equation~\eqref{eq:hg_force} is detailed in~\cite{Sun2022Understanding}. We also note that the sequence of element nodes differs slightly; readers should exercise caution when performing their own calculations.

\subsection*{$Q^2-Q^1$ case}
The shape functions of the $Q^1$ FEM space for the thermodynamic variables are listed in~\cite[Equation (5.1), Section 5]{Sun2025High}. The values of these four shape functions at the nine Gauss quadrature points form the interpolation matrix $M_{\text{interp}}$.
\begin{displaymath}
  \begin{bmatrix}
    0.7+1.5\sqrt{0.2} & -0.2 & -0.2 & 0.7-1.5\sqrt{0.2} \\
    0.25+0.75\sqrt{0.2} & 0.25+0.75\sqrt{0.2} & 0.25-0.75\sqrt{0.2} & 0.25-0.75\sqrt{0.2} \\
    -0.2 & 0.7+1.5\sqrt{0.2} & 0.7-1.5\sqrt{0.2} & -0.2\\
    0.25+0.75\sqrt{0.2} & 0.25-0.75\sqrt{0.2} & 0.25+0.75\sqrt{0.2} & 0.25-0.75\sqrt{0.2} \\
    0.25 & 0.25 & 0.25 & 0.25 \\
    0.25-0.75\sqrt{0.2} & 0.25+0.75\sqrt{0.2} & 0.25-0.75\sqrt{0.2} & 0.25+0.75\sqrt{0.2} \\
    -0.2 & 0.7-1.5\sqrt{0.2} &  0.7+1.5\sqrt{0.2} & -0.2 \\
    0.25-0.75\sqrt{0.2} & 0.25-0.75\sqrt{0.2} & 0.25+0.75\sqrt{0.2} & 0.25+0.75\sqrt{0.2} \\
    0.7-1.5\sqrt{0.2} & -0.2 & -0.2 & 0.7+1.5\sqrt{0.2}
  \end{bmatrix}  
\end{displaymath}
Let $\tilde{\rho}$ also denote the higher-order density variable, which belongs to the $Q^2$ FEM space. The values of $\tilde{\rho}$ can be updated by enforcing the mass conservation law at the nine Gauss-Legendre quadrature points, analogous to equation~\eqref{discrete_mass_conservation_Q1}.

The next task is to compute the term $\nabla N_j(\gv{\xi}_q)\det J(\gv{\xi}_q)$ in equation~\eqref{eq:hg_force}. While the momentum conservation law is discretized using a four-point Gauss-Legendre quadrature rule, the anti-hourglass force is constructed with a nine-point Gauss-Legendre quadrature. This means that $\nabla N_j(\gv{\xi}_q)\det J(\gv{\xi}_q)$ must be evaluated at different points. The advantage is that the special structure of the derivatives of the kinematic shape functions is preserved~\cite[Equation (4.14), Section 4]{Sun2025High}. As a result, these derivatives at the nine Gauss quadrature points can be reused in the computation of the artificial viscosity. 

For the $Q^3-Q^2$ case or higher-order cases, the calculation of the anti-hourglass force follows a similar procedure. This involves constructing the interpolation matrix from the shape functions of the thermodynamic variables, solving the mass conservation law at the higher-order quadrature points, and computing the derivatives of the kinematic shape functions. Finally, the anti-hourglass force is calculated using these auxiliary variables.

We conclude this section with a discussion on the robustness of the proposed framework. In our previous work, we pointed out that employing higher-order Gauss quadrature rules to solve the mass conservation law can compromise the robustness of the scheme. This observation also applies to the construction of the anti hourglass force. The higher-order density variable arises from the higher-order Jacobian determinants at the Gauss quadrature points, which are more prone to singularity due to grid distortion compared with the low-order case. However, this is a necessary trade-off; otherwise, the scheme would suffer from hourglass distortion.

An alternative approach is to construct the anti-hourglass force by limiting the velocity field, such as the viscosity-based hourglass control~\cite{Flanagan1981uniform} in the $Q^1-P^0$ case. In this method, the anti-hourglass force must be derived on a case by case basis. We found that such hourglass control is poorly compatible with the tensor viscosity in the $Q^1-P^0$ case~\cite{Sun2022On}. Future work will explore how to construct a suitable viscosity-based anti hourglass force for higher-order cases.
\section{The artificial viscosity}\label{sec:artificial_viscosity}
The artificial viscosity adopted in the proposed framework is consistent with that introduced in~\cite{Dobrev2012High}. Our study focuses on improving its computational efficiency and formulating a concise representation suitable for high-order schemes. Moreover, we show that the intermediate variable associated with the viscosity parameter and time-step control can be computed explicitly. Finally, we examine the contribution of artificial viscosity to the robustness of the proposed framework.

For the momentum equation, we consider only the artificial viscosity.
\begin{equation}\label{eq:viscosity_continuous}
\rho \frac{d \vec{u} }{d t} = \div (\mu \gv{\sigma}_{a}),
\end{equation}
$\mu$ denotes the viscosity parameter, which will be discussed in detail later. We present two forms of $\gv{\sigma}_{a}$. The first is the velocity-gradient formulation, corresponding to the classical expression of viscosity.
\begin{equation}\label{viscosity_form_asymmetry}
\gv{\sigma}_{a}= \nabla \vec{u} =
\begin{bmatrix}
   \frac{\partial u }{\partial r }  & \frac{\partial v }{\partial z }\\
   \frac{\partial v }{\partial r }  & \frac{\partial v }{\partial z },
\end{bmatrix}
\end{equation}
where $\vec{u} = [u, v]^{T}$. The second form is the symmetrized velocity gradient, as proposed in~\cite{Dobrev2012High}. 
\begin{equation}\label{viscosity_form_symmetry}
    \gv{\sigma}_{a}= \epsilon(\vec{u}) = \frac{1}{2} (\nabla \vec{u}+ \vec{u}\nabla)=
    \begin{bmatrix}
       \frac{\partial u }{\partial r }  & \frac{1}{2}(\frac{\partial v }{\partial r }+\frac{\partial v }{\partial z })\\
       \frac{1}{2}(\frac{\partial v }{\partial r }+\frac{\partial v }{\partial z })  & \frac{\partial v }{\partial z }
    \end{bmatrix}.
    \end{equation}
We first consider the classical expression of $\gv{\sigma}_{a}$ as defined in equation~\eqref{viscosity_form_asymmetry}. Accordingly, equation~\eqref{eq:viscosity_continuous} can be reformulated in explicit coordinate representation.
\begin{equation}\label{eq:viscosity_continuous_2}
    \begin{split}
    \rho \frac{d u }{d t} &= \nabla \cdot (\mu \nabla  u), \\
    \rho \frac{d v }{d t} &= \nabla \cdot (\mu \nabla  v).
    \end{split}
\end{equation}
For an element $K$, the discrete velocity components $u_h$ and $v_h$ are represented as linear combinations of the nodal values weighted by the associated shape functions.
\begin{displaymath}
\begin{split}
    u_h|_K&=\sum_{i=1}^{kdof}u_i N_i,\\
    v_h|_K&=\sum_{i=1}^{kdof}v_i N_i.
\end{split}    
\end{displaymath}

Applying the Petrov-Galerkin method, the weak form of equation~\eqref{eq:viscosity_continuous_2} is discretized as follows.
\begin{equation}\label{eq:viscosity_discretization}
    \begin{split}
        \sum_{i=1}^{kdof} \frac{d  u_{i}}{d t} \int_{K} \rho N_{i} N_j \dr\dz &= \int_{K} \mu \nabla u_h \cdot \nabla N_j \dr\dz, \\
        \sum_{i=1}^{kdof} \frac{d  v_{i}}{d t} \int_{K} \rho N_{i} N_j \dr\dz &= \int_{K} \mu \nabla v_h \cdot \nabla N_j \dr\dz
    \end{split}
\end{equation}
The left-hand side of equation~\eqref{eq:viscosity_discretization} represents the mass matrix, identical to that in the momentum conservation law. We focus on the right-hand side and, since the procedure is analogous in all coordinates, present only the $R$-direction computation.
\begin{equation}\label{eq:viscosity_rhs_discretization}
    \begin{split}
    \int_{K} \mu \nabla u_h \cdot \nabla N_j \dr\dz & = \int_{K} \mu \nabla \sum_{i=1}^{kdof}u_i N_i \cdot \nabla N_j \dr\dz \\
    &=\int_{K} \mu  \left ( \sum_{i=1}^{kdof}u_i \nabla  N_i \right) \cdot \nabla N_j \dr\dz \\
    &=\sum_{q=1}^{(m+1)^2}  \mu \omega_{q} \frac{1}{\det J(\gv{\xi}_q)} \left ( \sum_{i=1}^{kdof}u_i \nabla  N_i(\gv{\xi}_q)\det J(\gv{\xi}_q) \right) \cdot \nabla N_j(\gv{\xi}_q)\det J(\gv{\xi}_q)
    \end{split} 
\end{equation}
The form of equation~\eqref{eq:viscosity_rhs_discretization} differs slightly from that in~\cite{Kolev2009tensor,Sun2022On}, where the latter employs the matrix $J^{*T}J^{*}$, with $J^{*}$ denoting the adjoint of the Jacobian $J$. In the proposed framework, $\nabla (\gv{\xi}_q) \det J(\gv{\xi}_q)$ serves as a fundamental component that is used frequently. This choice preserves the compactness of both the equations and the code implementation, while maintaining a structure that allows the computation to be decomposed into a few simple terms. In our implementation, the values of $\nabla N_i(\gv{\xi}_q) \det J(\gv{\xi}_q)$ are first computed at the $(m+1)^2$ Gauss quadrature points, which can also be reused in the anti-hourglass force calculation. The gradient of the velocity field at the $(m+1)^2$ quadrature points can then be efficiently obtained by summation as follows:
\begin{displaymath}
    \nabla u_h (\gv{\xi}_q) =\sum_{i=1}^{kdof}u_i \left( \nabla  N_i(\gv{\xi}_q)\det J(\gv{\xi}_q) \right).
\end{displaymath}

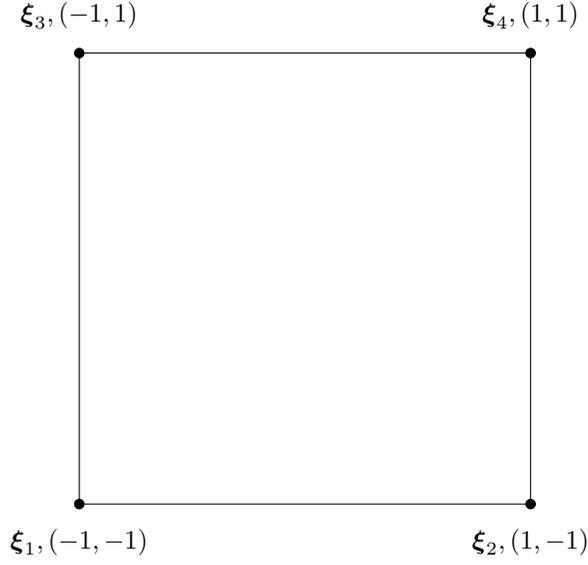
\begin{figure}
    \begin{center}
    \begin{tikzpicture}[scale = 1]
    \fill[black] (-3,-3) circle (2pt);
    \fill[black] (3,-3) circle (2pt);
    \fill[black] (3,3) circle (2pt);
    \fill[black] (-3,3) circle (2pt);
    
    \draw[black] (-3,-3) -- (3,-3);
    \draw[black] (3,-3) -- (3,3);
    \draw[black] (3,3) -- (-3,3);
    \draw[black] (-3,3) -- (-3,-3);
    \node[] at (-3,-3.5) {$\gv{\xi}_1,(-1,-1)$};
    \node[] at (3,-3.5)  {$\gv{\xi}_2,(1,-1)$};
    \node[] at (-3,3.5)  {$\gv{\xi}_3,(-1,1)$};
    \node[] at (3,3.5)   {$\gv{\xi}_4,(1,1)$};

    \end{tikzpicture}
    \caption{Reference element of the $Q^1$ FEM space, showing the degrees of freedom of the kinematic variables and the Gauss-Lobatto quadrature points (black circles), which are coincident.}
    \label{fig:Q_1_reference_element}
    \end{center}
\end{figure}

Next, we outline the derivation of the tensor viscosity~\cite{Campbell2001tensor,Kolev2009tensor,Sun2022On} within the proposed framework. Figure~\ref{fig:Q_1_reference_element} illustrates the Gauss-Lobatto quadrature points and the reference element of the $Q^1$ element, with the node sequence arranged according to the $R$-direction priority rule. Due to this ordering, the resulting form of the viscosity differs slightly from that presented in the original references. The following derivation details this formulation.

The gradients of the shape functions were derived in previous work~\cite{Sun2025High} [Equations (4.9)-(4.10), Section 4] and are here reformulated as follows:
 \begin{equation}\label{eq:Q1_velocity_shape_function_gradient_234}
     \begin{split}
        (\det J) \nabla N_1&
        =
    \begin{bmatrix}
        -N_1 z_{32}-N_2 z_{42}-N_3 z_{34}\\
        N_1 r_{32}+N_2 r_{42}+N_3 r_{34} 
    \end{bmatrix},\\
     (\det J) \nabla N_2&
          =
      \begin{bmatrix}
          -N_1 z_{13}+N_2 z_{41}-N_4 z_{34}\\
          N_1 r_{13}-N_2 r_{41}+N_4 r_{34} 
      \end{bmatrix},\\
      (\det J) \nabla N_3&
      =
     \begin{bmatrix}
      -N_1 z_{21}-N_3 z_{41}-N_4 z_{42}\\
      N_1 r_{21}+N_3 r_{41}+N_4 r_{42} 
     \end{bmatrix}, \\     
      (\det J) \nabla N_4&
      =
  \begin{bmatrix}
     -N_2 z_{21}-N_3 z_{13}+N_4 z_{32}\\
     N_2 r_{21}+N_3 r_{13}-N_4 r_{32}
  \end{bmatrix}.
     \end{split}
  \end{equation}

For the four Gauss-Lobatto quadrature points,
\begin{displaymath}
N_i(\gv{\xi_q})=\delta_{iq}, \quad \delta_{iq}=1 , \ i =q, \ \delta_{iq}=0, \ i\neq q.
\end{displaymath}

Subsequently, the gradient of $N_1 \det J$ may be explicitly represented as follows:
\begin{equation}
\nabla N_1 (\gv{\xi_1}) \det J(\gv{\xi_1})  =
\begin{bmatrix}
    -z_{32}\\
   r_{32} 
\end{bmatrix},
\nabla N_1 (\gv{\xi_2}) \det J(\gv{\xi_2}) =
\begin{bmatrix}
    -z_{42}\\
   r_{42} 
\end{bmatrix},
\nabla N_1 (\gv{\xi_3}) \det J(\gv{\xi_3}) =
\begin{bmatrix}
    -z_{34}\\
   r_{34} 
\end{bmatrix}.
\end{equation}
where $\nabla N_1 (\gv{\xi_4}) \det J(\gv{\xi_4}) = [0, 0]^T$, and the gradients of the remaining shape functions can be derived analogously.

The gradient of $u_h$, denoted $\nabla u_h$, can be derived as follows:,
\begin{equation}
\begin{split}
\nabla u_h(\gv{\xi}_1) \det J(\gv{\xi}_1) &=
u_1
\begin{bmatrix}
    -z_{32}\\
   r_{32} 
\end{bmatrix}
+
u_2
\begin{bmatrix}
    -z_{13}\\
   r_{13} 
\end{bmatrix}
+
u_3
\begin{bmatrix}
    -z_{21}\\
   r_{21} 
\end{bmatrix}, \\
\nabla u_h(\gv{\xi}_2) \det J(\gv{\xi}_2)  &=
u_1
\begin{bmatrix}
    -z_{42}\\
   r_{42} 
\end{bmatrix}
+
u_2
\begin{bmatrix}
    z_{41}\\
   -r_{41} 
\end{bmatrix}
+
u_4
\begin{bmatrix}
    -z_{21}\\
   r_{21} 
\end{bmatrix}, \\
\nabla u_h(\gv{\xi}_3) \det J(\gv{\xi}_3)  &=
u_1
\begin{bmatrix}
    -z_{34}\\
   r_{34} 
\end{bmatrix}
+
u_3
\begin{bmatrix}
    -z_{41}\\
   r_{41} 
\end{bmatrix}
+
u_4
\begin{bmatrix}
    -z_{13}\\
   r_{13} 
\end{bmatrix}, \\
\nabla u_h(\gv{\xi}_4) \det J(\gv{\xi}_4)  &=
u_2
\begin{bmatrix}
    -z_{34}\\
   r_{34} 
\end{bmatrix}
+
u_3
\begin{bmatrix}
    -z_{42}\\
   r_{42} 
\end{bmatrix}
+
u_4
\begin{bmatrix}
    z_{32}\\
   -r_{32} 
\end{bmatrix}, 
\end{split}
\end{equation}

Finally, the nodal force in the $R$-direction for node $1$ is given by:
\begin{equation}\label{Q_1_viscosity_force_1R}
    \begin{split}
        \int_{K} \mu \nabla u_h \cdot \nabla N_1 \dr\dz &= \sum_{q=1}^{4} \mu \frac{1}{\det J(\gv{\xi}_q)} \nabla u_h(\gv{\xi}_q)\det J(\gv{\xi}_q) \cdot \nabla N_1(\gv{\xi}_q) \det J(\gv{\xi}_q) \\
        & =\mu_1 \frac{1}{\det J (\gv{\xi}_1)} \left[  u_1(r_{32}^2+z_{32}^2) +u_2(r_{13} r_{32} + z_{13} z_{32}) +u_3 (r_{21}r_{32}+z_{21} z_{32}) \right] \\
        &+\mu_2 \frac{1}{\det J (\gv{\xi}_2)} \left[ u_1(r_{42}^2+z_{42}^2) +u_2 (-r_{41}r_{42}-z_{41}z_{42}) +u_4(r_{21}r_{42}+z_{21}z_{42})     \right] \\
        &+ \mu_3 \frac{1}{\det J (\gv{\xi}_3)} \left[u_1(r_{34}^2+z_{34}^2) +u_3(r_{41}r_{34}+z_{41}z_{34})+u_4(r_{13}r_{34}+z_{13}z_{34})     \right]. 
    \end{split}
\end{equation}
Equation~\eqref{Q_1_viscosity_force_1R} provides an explicit expression for the viscosity force, as given in~\cite{Sun2022On} [Equation (2.19), Section 2]. Other forms of tensor viscosity can be obtained by employing different quadrature rules. The form~\eqref{Q_1_viscosity_force_1R} is explicit and many variables can be reused, but it is not suitable for higher-order cases. In practical implementations, such formula requires too many temporary variables, which is not conducive to code readability and maintenance. Therefore, we recommend using the more concise form presented in equation~\eqref{eq:viscosity_rhs_discretization} with precomputed $\nabla u_h,\nabla v_h$. Despite the slight increase in computational cost, this approach significantly enhances code clarity and maintainability.

For the general $Q^m-Q^{m-1}$ FEM space pair, the implementation of artificial viscosity follows a similar procedure. We employ $(m+1)^2$-point Gauss-Legendre quadrature rules to evaluate the right-hand side of equation~\eqref{eq:viscosity_discretization}. The derivatives of the shape functions, $\nabla N_i \det J$, are computed at the $(m+1)^2$ quadrature points, preserving a specified structure. The gradient of the velocity field is then obtained by summation. In the symmetrized case, $\gv{\sigma}_{a} = \epsilon(\vec{u})$, the discretization scheme proceeds analogously. The remaining question is how to determine the viscosity parameter.

\subsection{The artificial viscosity parameter}
From the preceding calculations, it follows that the artificial viscosity parameter $\mu$ is defined at each quadrature point, specifically at the $(m+1)^2$ Gauss-Legendre quadrature points. The general form of $\mu$ is given by:
\begin{equation}\label{viscosity_parameter}
\mu=\rho(c_1 c_{s} l_{c}  +c_2 l_{c}^2 |\Delta \vec{u}| ).    
\end{equation}
Here, $c_1$ and $c_2$ are non-dimensional, nonnegative linear and quadratic scaling parameters, respectively; $c_s$ denotes the sound speed, $\rho$ the density, and $l_c$ the characteristic length. The term $|\Delta \vec{u}|$ serves as an indicator of compression. One limitation of the viscosity parameter $\mu$ in equation~\eqref{viscosity_parameter} is that the limiter is neglected.

The compression indicator $|\Delta \vec{u}|$ generally determines whether the viscosity parameter $\mu$ is active or vanishes.
\begin{displaymath}
    \mu=0, \quad  \text{if} \ |\Delta \vec{u}|\ge 0.
\end{displaymath}
For the $Q^1-P^0$ case, the divergence of the velocity field constitutes a suitable candidate.
\begin{displaymath}
    |\Delta \vec{u}|=\div \vec{u}.
\end{displaymath}
For higher-order FEM space pairs, however, the numerical performance of this indicator is unsatisfactory. Following the approach in~\cite{Dobrev2012High}, we adopt the minimum eigenvalue of the symmetrized gradient, $\epsilon(\vec{u})$, as the indicator. In the two-dimensional case, both this eigenvalue and its corresponding eigenvector can be computed analytically. Consider the following eigenvalue problem:
\begin{displaymath}
    \lambda I - \epsilon(\vec{u})=
    \begin{bmatrix}
        \lambda - \frac{\partial u }{\partial r }  & -\frac{1}{2}(\frac{\partial v }{\partial r }+\frac{\partial v }{\partial z })\\
        -\frac{1}{2}(\frac{\partial v }{\partial r }+\frac{\partial v }{\partial z })  & \lambda-\frac{\partial v }{\partial z }
     \end{bmatrix}.
\end{displaymath}
The minimum eigenvalue $\lambda_{m}$ is
\begin{equation}\label{eigenvalue}
    \lambda_{m}= \frac{1}{2}(\frac{\partial u }{\partial r } + \frac{\partial v }{\partial z }) + \frac{1}{2}\sqrt{(\frac{\partial u }{\partial r } - \frac{\partial v }{\partial z })^2 +(\frac{\partial u }{\partial z }+\frac{\partial v}{\partial r })^2}.
\end{equation}
The corresponding eigenvector $\vec{e}$ is
\begin{equation}
    \vec{e}=
    \begin{bmatrix}
        \frac{1}{2}(\frac{\partial v }{\partial r }+\frac{\partial v }{\partial z }) \\
        \lambda_m - \frac{\partial u }{\partial r }     
\end{bmatrix} 
    =
    \begin{bmatrix}
        \frac{1}{2}(\frac{\partial v }{\partial r }+\frac{\partial v }{\partial z }) \\
        \frac{1}{2}(-\frac{\partial u }{\partial r } + \frac{\partial v }{\partial z }) + \frac{1}{2}\sqrt{(\frac{\partial u }{\partial r } - \frac{\partial v }{\partial z })^2 +(\frac{\partial u }{\partial z }+\frac{\partial v}{\partial r })^2}
    \end{bmatrix}   
\end{equation}

The eigenvalue $\lambda_m$ is taken as the indicator of compression.
\begin{equation}\label{eigenvalue_indicator}
    |\Delta \vec{u}|=\lambda_{m}.
\end{equation}
In the preceding viscosity computation, once $\nabla u_h$ and $\nabla v_h$ are obtained, the eigenvalue $\lambda_m$ and the corresponding eigenvector $\vec{e}$ can be computed directly. 

The computation of the characteristic length $l_c$ is determined by the eigenvector $\vec{e}$ and the Jacobian matrix.
\begin{equation}\label{characteristic_length}
l_c=l^0 \frac{|(J^0)^{-1}J \vec{e}|}{|\vec{e}|},    
\end{equation}
$l^0$ denotes the initial characteristic length at $t = 0$, and $J^0$ represents the Jacobian matrix at $t = 0$.

Assume that $J^0$ and $J$ have the following elements:
\begin{displaymath}
J^0=
\begin{bmatrix}
    J^0_{11} & J^0_{12} \\
    J^0_{21} & J^0_{22}
\end{bmatrix},  
J=
\begin{bmatrix}
    J_{11} & J_{12} \\
    J_{21} & J_{22}
\end{bmatrix}.,  
\end{displaymath}
Then, the inverse matrix $(J^0)^{-1}$ is given by:
\begin{displaymath}
    (J^0)^{-1}=\frac{1}{\det J^0}
    \begin{bmatrix}
        J^0_{22} & -J^0_{12} \\
        -J^0_{21} & J^0_{11}
    \end{bmatrix},  
\end{displaymath}
The matrix $(J^0)^{-1} J$ is expressed as:
\begin{displaymath}
    (J^0)^{-1} J=\frac{1}{\det J^0}
    \begin{bmatrix}
        J^0_{22}J_{11}-J^0_{12}J_{21} & J^0_{22}J_{12}-J^0_{12}J_{22} \\
        J^0_{11}J_{21}-J^0_{21}J_{11} & J^0_{11}J_{22}-J^0_{21}J_{12}
    \end{bmatrix}=\frac{1}{\det J^0}
    \begin{bmatrix}
        M_{11} & M_{12}\\
        M_{21} & M_{22}
    \end{bmatrix},     
\end{displaymath}

For $l^0$, we approximate it using the square root of the initial Jacobian determinant, $\sqrt{\det J^0}$.
\begin{equation}
    l^0=\sqrt{\det J^0}.
\end{equation}

The eigenvector is normalized and denoted by $\vec{n}_e$.
\begin{displaymath}
    \vec{n}_e=\frac{\vec{e}}{|e|}=
    \begin{bmatrix}
        ne_{1}\\
        ne_{2}
    \end{bmatrix}.
\end{displaymath}
The characteristic length $l_c$ at each quadrature point $\gv{\xi}_q$ can then be computed using the following formula:
\begin{equation}
    l_c=\sqrt{\frac{\omega_q}{\det J^0}[M_{11}^2+M_{12}^2+M_{21}^2+M_{22}^2+2*(M_{11}M_{12}+M_{21}M_{22})ne_{1}ne_{2}]}.
\end{equation}

In~\cite{Dobrev2012High}, a parameter $c_{\rm vor}$, related to vorticity, is introduced. It is then multiplied with the linear term to obtain the complete expression for $\mu$.
\begin{equation}\label{viscosity_vorticity_parameter}
    \mu=\rho(c_1 c_{vor} c_{s} l_{c}  +c_2 l_{c}^2 |\Delta \vec{u}| ),    
\end{equation}
where 
\begin{displaymath}
c_{vor}= \frac{|\div \vec{u}|}{ \|\grad \vec{u}\|}.    
\end{displaymath}
In two dimensional case,
\begin{displaymath}
    c_{vor}=\frac{|\frac{\partial u }{\partial r } + \frac{\partial v }{\partial z }|}{\sqrt{\frac{\partial u }{\partial r }^2+ \frac{\partial u }{\partial z } ^2 +\frac{\partial v }{\partial r }^2 +\frac{\partial v }{\partial z }^2 }}.
\end{displaymath}
Such derivatives as $\frac{\partial u }{\partial r }$, is calculated previously.

For the thermodynamic variables, such as $c_s$, the values are defined at $m^2$ quadrature points. These values need to be interpolated from the $m^2$ points to the $(m+1)^2$ points. The interpolation matrix $M_{interp}$, defined in Section~\ref{sec:hourglass}, is reused for this purpose. For the density variable $\rho$, the existing term $\tilde{\rho}$, also defined in Section~\ref{sec:hourglass}, serves as a suitable candidate.

\subsection{Time step control}
The time step control must also be computed at each of the $(m+1)^2$ quadrature points within each element. It is generally determined using the following formula:
\begin{displaymath}
\left(\frac{c_s+Q}{l_{\tau}}\right)^{-1},
\end{displaymath}
$c_s$ denotes the previously defined sound speed, $l_{\tau}$ is the characteristic length used for time step control, and $Q$ is a viscosity related term, which typically takes the following form:
\begin{displaymath}
    Q=\frac{\mu}{\rho l_{\tau}},
\end{displaymath}
$\mu$ denotes the previously defined viscosity parameter. The above formula can then be simplified as follows:
\begin{displaymath}
    \left(\frac{c_s}{l_{\tau}}+\frac{\mu}{\rho l_{\tau}^2}\right)^{-1}   
\end{displaymath}

Following the approach in~\cite{Dobrev2012High}, we take the minimum singular value of the Jacobian matrix as the characteristic length $l_{\tau}$. This definition can be expressed explicitly. In two dimensions, the Jacobian matrix is given by:
\begin{displaymath}
J=
\begin{bmatrix}
    J_{11} & J_{12} \\
    J_{21} & J_{22}
\end{bmatrix},
\end{displaymath}
Then the square of minimum singular value of $J$ is
\begin{displaymath}
    l_{\tau}^2=\frac{1}{2}(J_{11}^2+J_{12}^2+J_{21}^2+J_{22}^2)-\frac{1}{2}\sqrt{(J_{11}^2+J_{21}^2-J_{12}^2-J_{22}^2)^2+4(J_{11}J_{12}+J_{21}J_{22})^2}.
\end{displaymath}

The time-step candidate is taken as the minimum value across all elements and quadrature points.
\begin{displaymath}
    \tau=\text{CFL} \min_{\gv{\xi}_q}  \left(\frac{c_s}{l_{\tau}}+\frac{\mu}{\rho l_{\tau}^2}\right)^{-1}.
\end{displaymath}

\subsection{The robustness of the artificial viscosity}
At the end of this section, we discuss the robustness of the artificial viscosity. In principle, one is free to choose any quadrature rule to solve the viscosity equation~\eqref{eq:viscosity_continuous}. However, the discretization of such artificial viscosity inevitably introduces instabilities. Recall the discretized formula of the artificial viscosity in equation~\eqref{eq:viscosity_rhs_discretization}, which contains a term involving the reciprocal of the Jacobian determinant.
\begin{displaymath}
    \frac{1}{\det J (\gv{\xi}_q)}
\end{displaymath}
If the Jacobian determinant at a quadrature point approaches zero, the artificial viscosity force becomes unbounded, which is the source of instability. It should be noted that the Jacobian determinant also affects the anti-hourglass force. A singular Jacobian can lead to an infinite density variable, resulting in incorrect pressure variations and, consequently, an inappropriate computation of the anti-hourglass force. As discussed above, the Jacobian determinant depends on both the quadrature point location and the element shape. Higher-order quadrature rules are more sensitive to grid distortion, exacerbating these issues.

For the $Q^{m}-Q^{m-1}$ FEM space pair, the $m^2$- or $(m+1)^2$-point Gauss-Legendre quadrature rules are optimal choices for discretizing the artificial viscosity. Using Gauss-Legendre quadrature rules with precision lower than $m^2$ points sacrifices the accuracy of the velocity field unnecessarily, since the mass conservation law is already solved using the $m^2$-point rule, whose robustness is also constrained by the Jacobian determinant. Conversely, quadrature rules of order higher than $(m+1)^2$ points are unnecessary and may even compromise robustness. The $(m+1)^2$-point rule provides sufficient precision for accurately computing the artificial viscosity.

The choice between these two quadrature rules is subtle, due to the interplay between artificial viscosity and hourglass control. Using the $m^2$-point rule sacrifices efficiency in order to gain robustness. However, the artificial viscosity is then computed less accurately, reducing its effectiveness in suppressing hourglass motion. Conversely, if the anti-hourglass force is computed with the $(m+1)^2$-point rule, robustness may also be limited. These quadrature rules may be better coordinated with hourglass control methods based on velocity limiting, which will be explored in future work.

In the proposed framework, the $(m+1)^2$-point Gauss-Legendre quadrature rule is preferred, as it computes the artificial viscosity accurately, thereby aiding in the suppression of hourglass motion. For the $Q^1-P^0$ case, the principles governing the efficiency and robustness of the artificial viscosity were discussed in previous work~\cite{Sun2022On}. Both artificial viscosity and hourglass control are subject to the same instabilities, highlighting the need to balance efficiency and robustness.

\section{Energy conservation for hourglass and artificial viscosity}\label{sec:energy_conservation}
In the previous section, we introduced the computation of the anti-hourglass force and the artificial viscosity force for the momentum equation. For energy conservation, however, the work done by both the anti-hourglass force and the artificial viscosity force must also be incorporated into the energy equation.

The continuous energy equation for the hourglass control and for the artificial viscosity are considered separately.
\begin{equation}\label{eq:energy_hg_av}
    \begin{split}
        \rho \frac{d e }{d t} &= \delta p \div{\vec{u}}\\
        \rho \frac{d e }{d t} &=\mu \gv{\sigma}_{a} : \grad{\vec{u}}
    \end{split}
\end{equation}
$\delta p$ denotes the pressure variation in equation~\eqref{pressure_variation}, and the symbol $:$ represents the summation convention.

As introduced in previous work~\cite{Sun2025High}, the discretization of equation~\eqref{eq:energy_hg_av} is more efficient when using the forces computed from the momentum equation discretization. In the previous section, we calculated the anti-hourglass force and the artificial viscosity force at the $(m+1)^2$ quadrature points. Let $\vec{f}_h^{,j}$ and $\vec{f}_v^{,j}$ denote the force vectors of the anti-hourglass and viscosity forces at each quadrature point, respectively, where the subscript $j$ corresponds to the degrees of freedom of the kinematic variable. The anti-hourglass and artificial viscosity forces at each quadrature point can then be expressed as follows:
\begin{equation}\label{eq:hg_force_denotation}
\vec{fh}_{j}(\gv{\xi}_q)= \delta p (\gv{\xi}_q) \nabla N_j (\gv{\xi}_q) \det J (\gv{\xi}_q),
\end{equation}  

\begin{equation}\label{eq:vis_force}
    \vec{fv}_{j}(\gv{\xi}_q)=
    \begin{bmatrix}
        \mu  \frac{1}{\det J(\gv{\xi}_q)} \left ( \sum_{i=1}^{kdof}u_i \nabla  N_i(\gv{\xi}_q)\det J(\gv{\xi}_q) \right) \cdot \nabla N_j(\gv{\xi}_q)\det J(\gv{\xi}_q)\\
        \mu \frac{1}{\det J(\gv{\xi}_q)} \left ( \sum_{i=1}^{kdof}v_i \nabla  N_i(\gv{\xi}_q)\det J(\gv{\xi}_q) \right) \cdot \nabla N_j(\gv{\xi}_q)\det J(\gv{\xi}_q)
    \end{bmatrix}.
\end{equation}

Since both forces employ the same quadrature rule, we combine them and denote the resulting force by $\vec{f}_{j}$.
\begin{equation}
    \vec{f}_{j}(\gv{\xi}_q)=\vec{fh}_{j}(\gv{\xi}_q)+\vec{fv}_{j}(\gv{\xi}_q).
\end{equation}

The discretization of the energy equation, considering only the anti-hourglass and viscosity forces, is then given by:
\begin{equation}\label{eq:energy_hg_av_discrete}
\sum_{k=1}^{tdof} \frac{d e_k}{d t}\int_{K} \rho_h  \phi_k \phi_l \dr\dz  =  \sum_{q=1}^{(m+1)^2} \omega_{q} \left(\sum_{j=1}^{kdof}\vec{u}_j \cdot \vec{f}_{j} (\gv{\xi}_q)\right)  \phi_l(\gv{\xi}_q)
\end{equation}
where $\phi_{l}$ denotes the shape function associated with the thermodynamic variables.

It is straightforward to verify that equation~\eqref{eq:energy_hg_av_discrete} represents the discretization of equation~\eqref{eq:energy_hg_av}. The $(m+1)^2$-point Gauss-Legendre quadrature rule is employed to evaluate the right-hand side. The notation $\vec{f}_{j}$ is introduced because the force vector is already computed in the momentum equation, eliminating the need for recalculation in the energy equation.


\section{Numerical example}\label{sec:numerical_results}

\subsection*{2D Taylor-Green vortex}
\begin{figure}
    \begin{center}
      \includegraphics[width=0.49\textwidth]{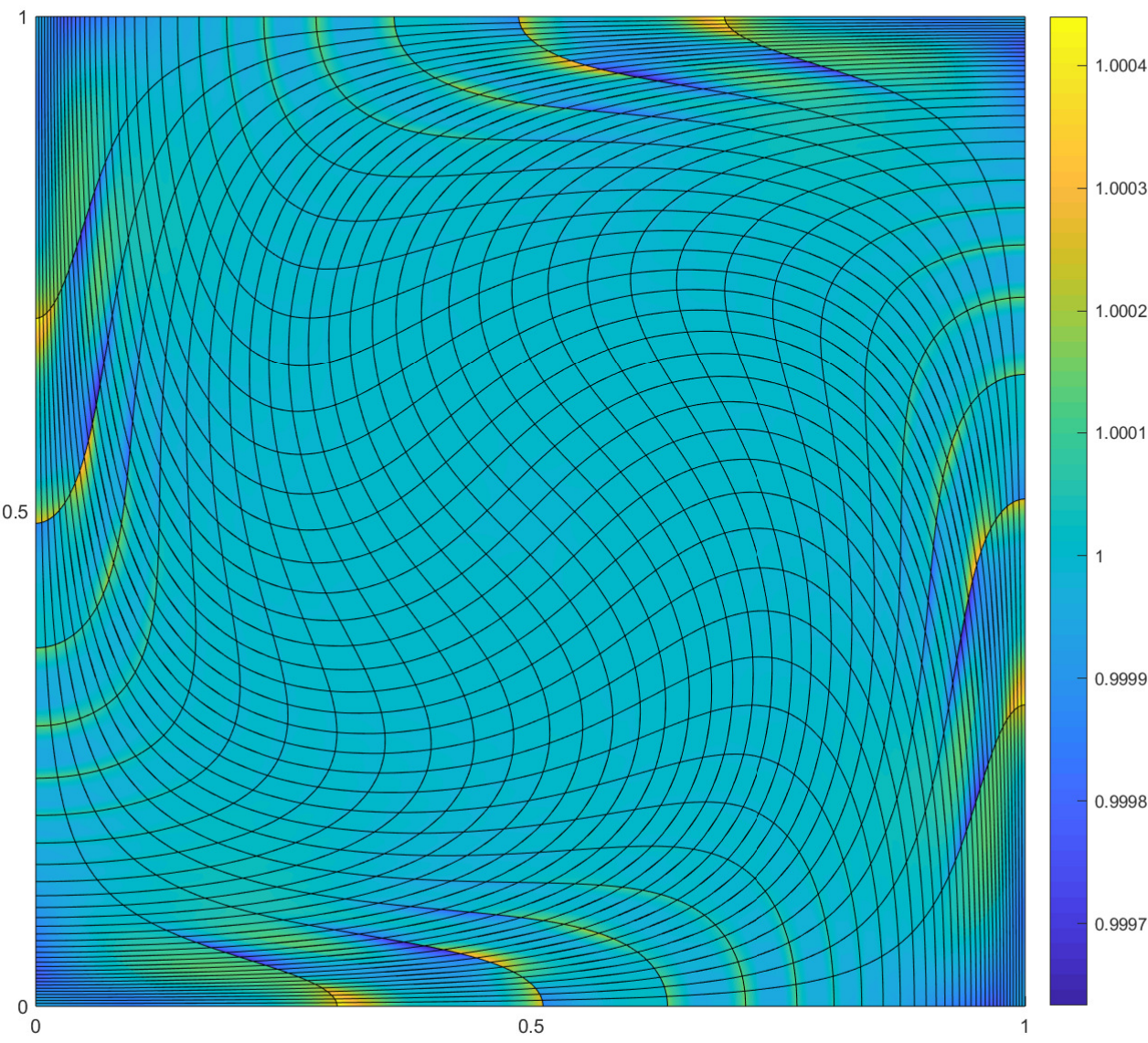}
      \includegraphics[width=0.49\textwidth]{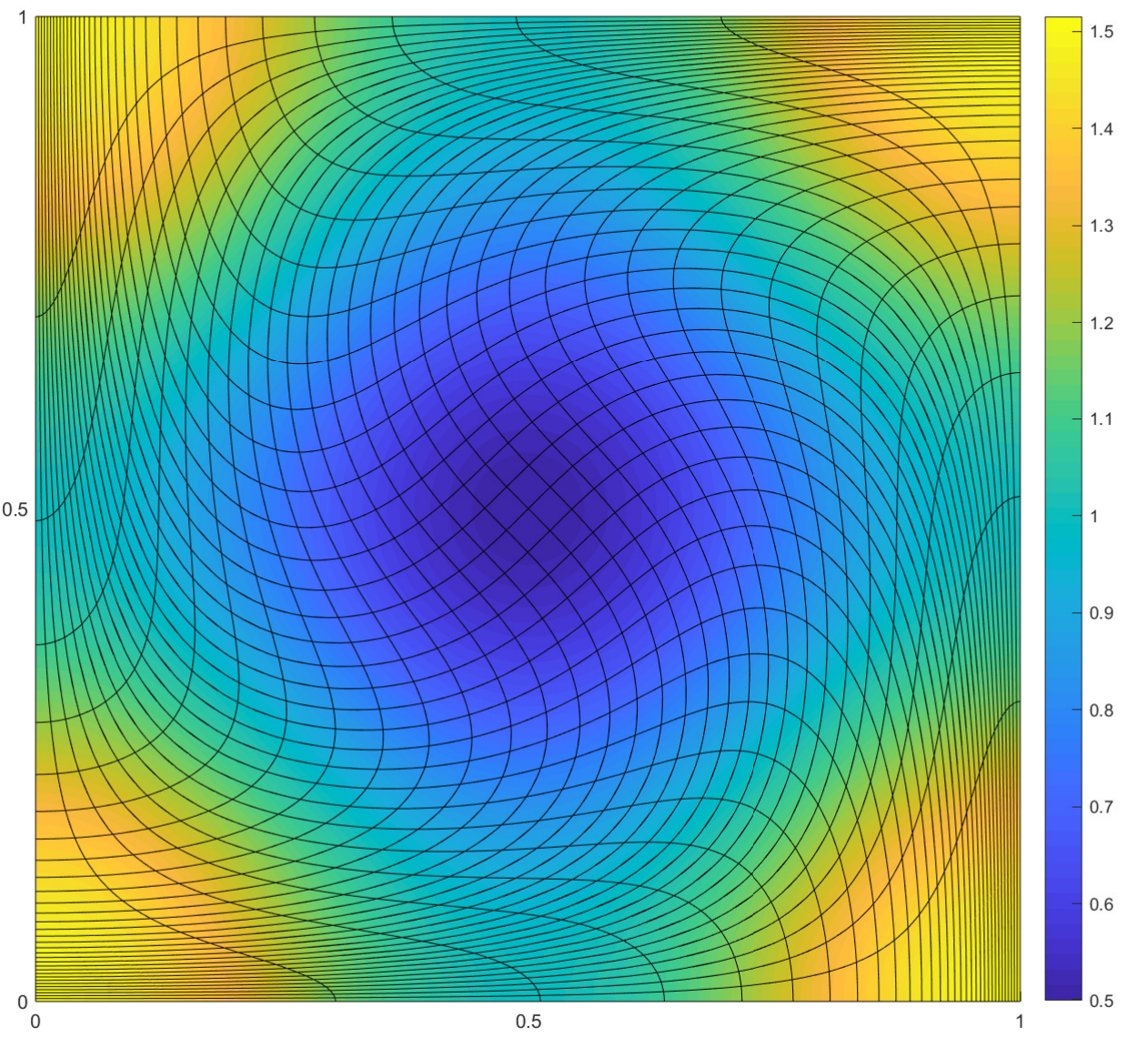}
      \caption{Numerical results of the Taylor-Green vortex using the $Q^3-Q^2$ element pair with $h = 1/32$ at $t = 0.75$: (a) density field, (b) pressure field.}
      \label{Taylor_Green_vortex_field}
    \end{center}
  \end{figure}

\begin{figure}
    \begin{center}
      \includegraphics[width=0.49\textwidth]{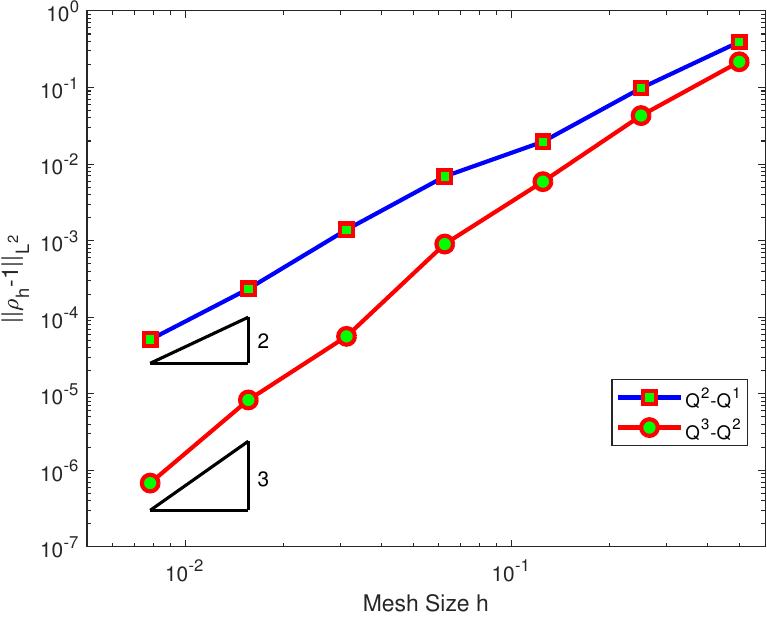}
      \includegraphics[width=0.49\textwidth]{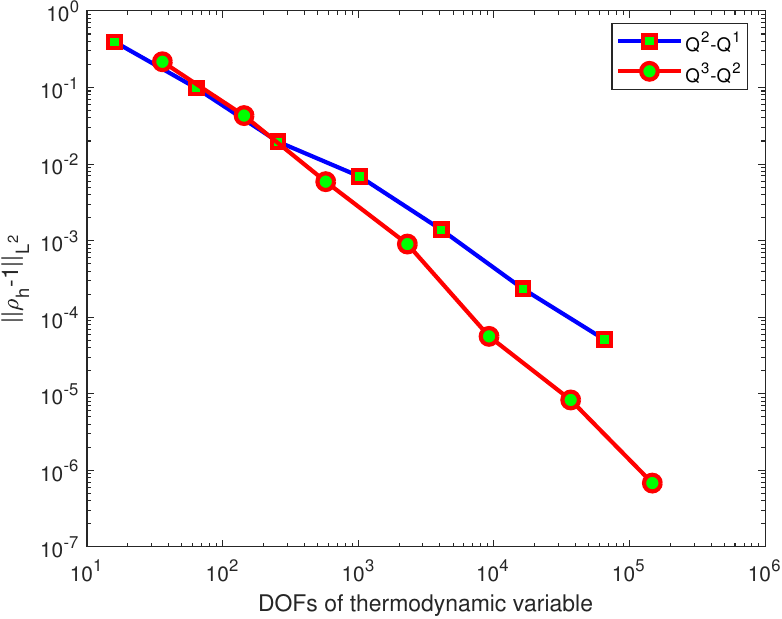}
      \caption{Numerical results of the Taylor-Green vortex at $t = 0.75$: (a) error versus mesh size $h$, (b) error versus degrees of freedom of the thermodynamic variable.}
      \label{Taylor_Green_vortex_error}
    \end{center}
  \end{figure}

The 2D Taylor-Green vortex is a smooth, shock-free problem. In previous work~\cite{Sun2025High}, results without hourglass control were presented at $t = 0.5$. Here, we show results with hourglass control but without artificial viscosity at $t = 0.75$, with the problem setup identical to~\cite{Sun2025High}. The error norm is the $L^2$ norm of the density variation. Figure~\ref{Taylor_Green_vortex_field} displays the density and pressure fields, while Figure~\ref{Taylor_Green_vortex_error}(a) and Table~\ref{tab:density_error_t.75} report the error norm and convergence order, which closely match theoretical predictions. The hourglass control preserves numerical accuracy effectively.

Figure~\ref{Taylor_Green_vortex_error}(b) highlights the advantage of high-order methods, showing that higher-order schemes achieve superior accuracy even with fewer degrees of freedom (DOFs). For $h = 1/64$ and the $Q^2-Q^1$ space pair, the DOFs for kinematic and thermodynamic variables are 16,641 and 16,384, respectively, yielding higher errors than the $Q^3-Q^2$ pair with $h = 1/32$ (9,409 and 9,216 DOFs). The $Q^3-Q^2$ results at $h = 1/32$ are comparable to the $Q^2-Q^1$ results at $h = 1/128$, despite the latter having roughly seven times more DOFs.

\begin{table}[]
    \centering
    \caption{Convergence rates of the density variation for the Taylor-Green vortex at $t = 0.75$.}
    \label{tab:density_error_t.75}
    \scalebox{1.0}{
    \begin{tabular}{|c|c|c|c|c|c|c|}
\hline 
\multirow{2}{*}{$h$} & \multicolumn{3}{c|}{$Q^2-Q^1$} &  \multicolumn{3}{c|}{$Q^3-Q^2$} \\
\cline{2-7}
      &DOFs & Error & Order&DOFs&  Error & Order\\
\hline
$\frac{1}{2}$ &(25,16) &3.9465E-1 & - &(49,36)& 2.1684E-1 & - \\
 \hline
$\frac{1}{4}$ &(81,64) &9.8200E-2 &2.0064 &(169,144) &4.2974E-2 &2.3351\\
 \hline
$\frac{1}{8}$ &(289,256) &1.9592E-2 &2.3258 &(625,576) &5.8916E-3 &2.8667\\
 \hline
$\frac{1}{16}$ &(1089,1024) &6.8888E-3 &1.5079 &(2401,2304) &9.0218E-4 &2.7071\\
 \hline
$\frac{1}{32}$ &(4225,4096) &1.3901E-3 &2.3090 &(9409,9216) &5.6164E-5 &4.0056\\
 \hline
 $\frac{1}{64}$ &(16641,16384) &2.3463E-4 &2.5667 &(37249,36864) &8.2973E-6 & 2.7589\\
 \hline 
 $\frac{1}{128}$ &(66049,65536) &5.0908e-05 &2.2044 &(148225,147456) &6.8112e-07& 3.60664\\
 \hline 
    \end{tabular}}
  \end{table}

\subsection*{Noh} 
An ideal gas moves toward the origin with unit velocity, generating a strong shock at the origin that propagates outward. The initial computational domain is $[0,1]^2$, with the equation-of-state parameter $\gamma = 5/3$, initial density $\rho = 1$, and initial internal energy $1 \times 10^{-10}$. This problem has an analytical solution: the shock reaches a radius of $0.2$ at $t = 0.6$, with post-shock density and pressure of $16$ and $5.33$, respectively. The mesh size is set to $h = 0.05$, and both hourglass control and artificial viscosity are enabled.

\begin{figure}
  \begin{center}
    \includegraphics[width=0.49\textwidth]{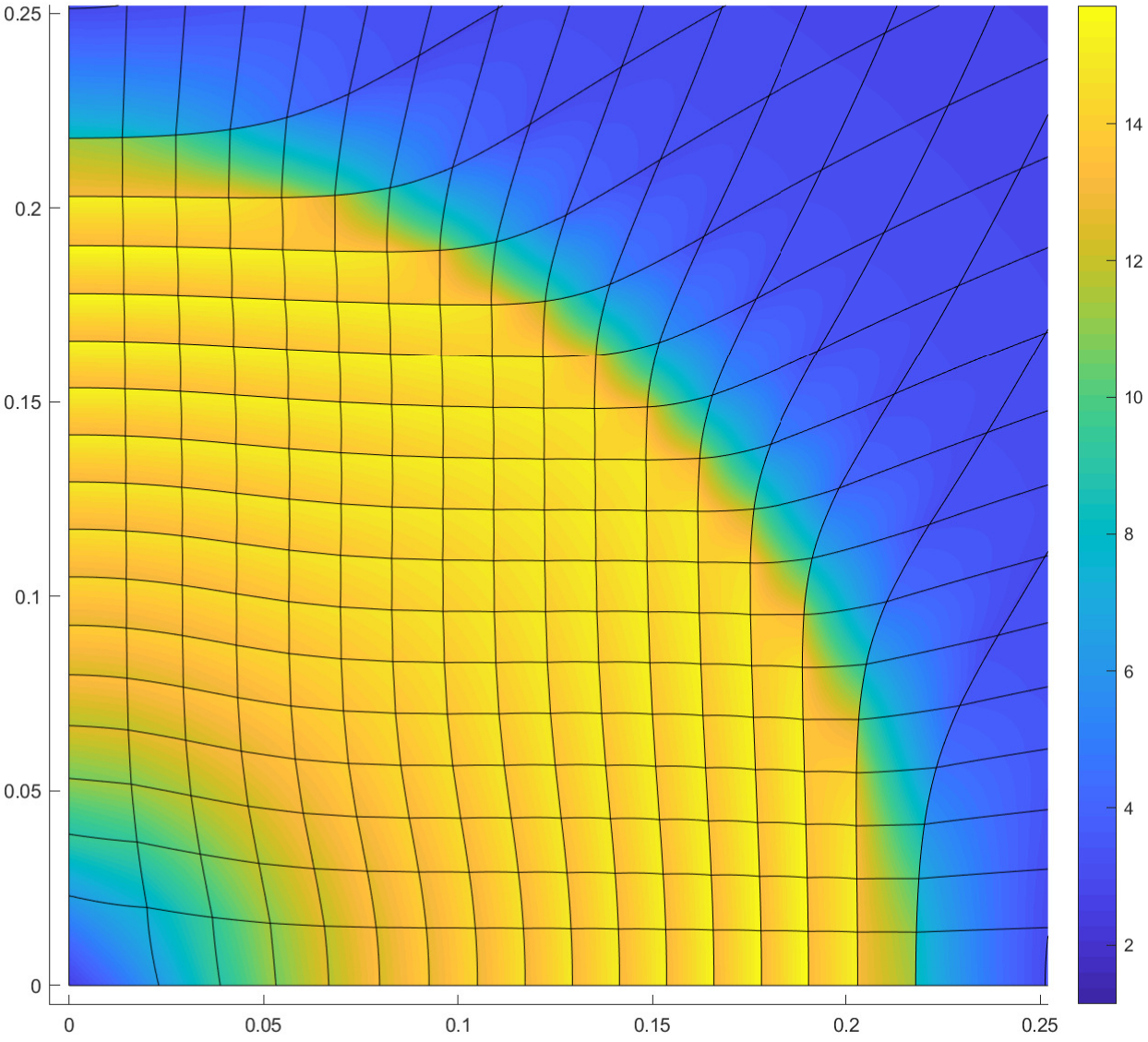}
    \includegraphics[width=0.49\textwidth]{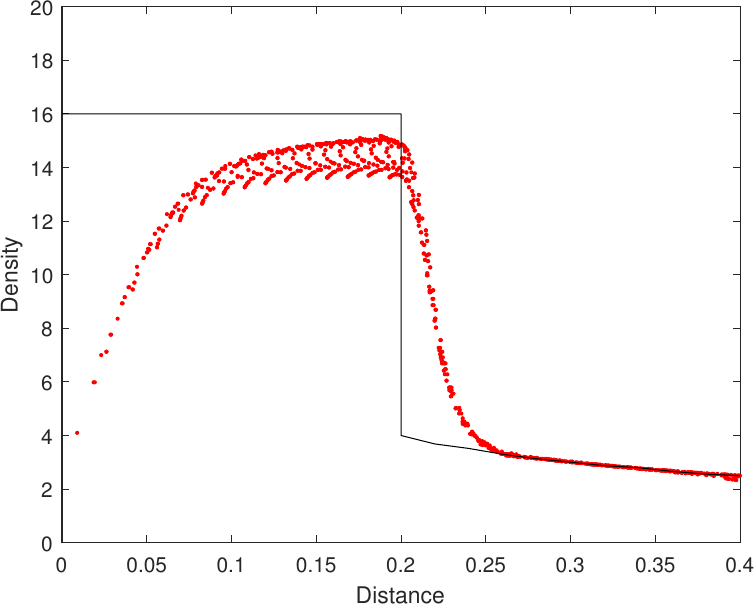}
    \caption{Numerical results of the Noh problem using the $Q^2-Q^1$ element pair with $h = \frac{1}{20}$ at $t = 0.6$: (a) density field, (b) density values compared with the analytical solution.}
    \label{Noh_Q2_Q1}
  \end{center}
\end{figure}

\begin{figure}
  \begin{center}
    \includegraphics[width=0.49\textwidth]{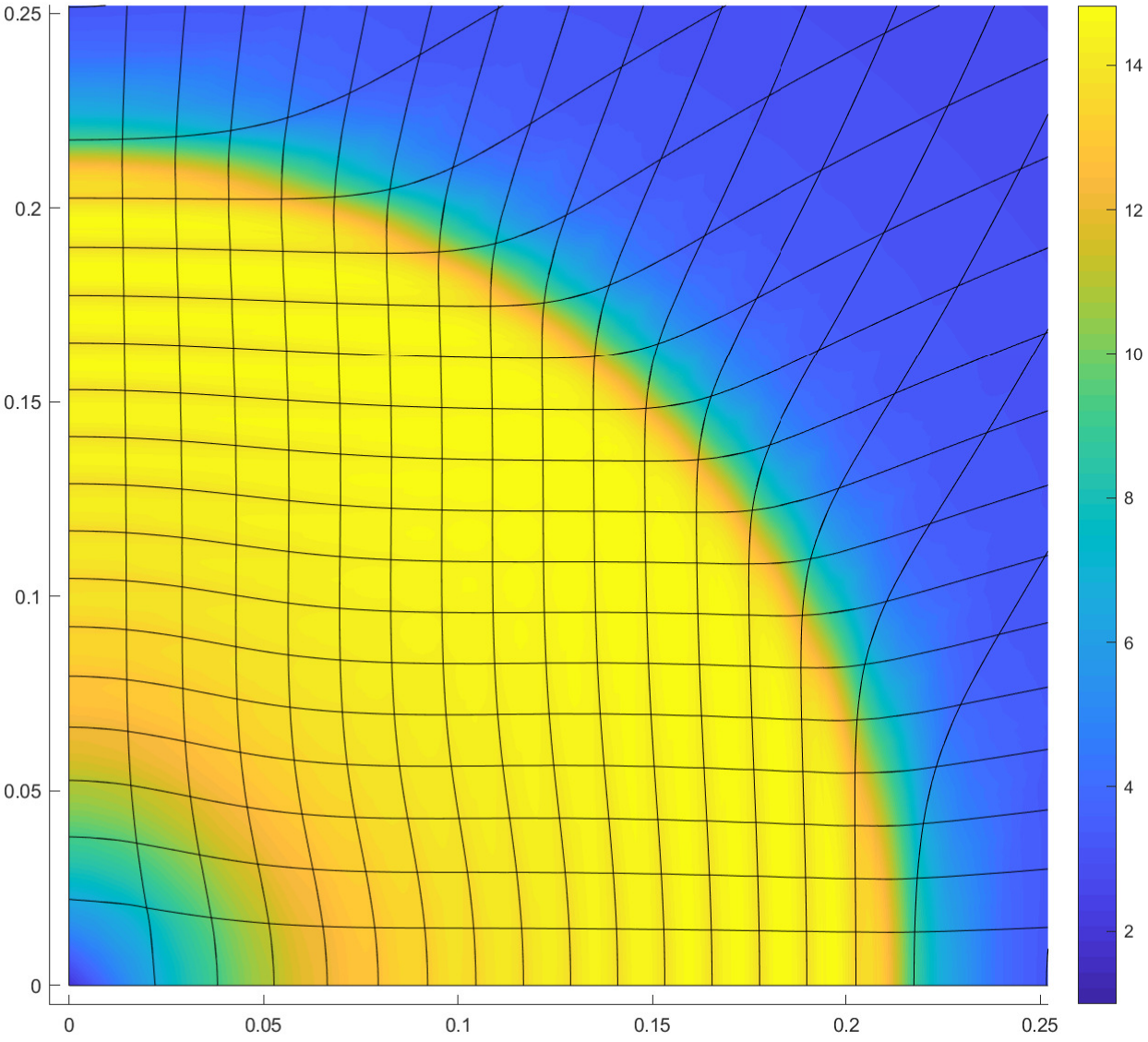}
    \includegraphics[width=0.49\textwidth]{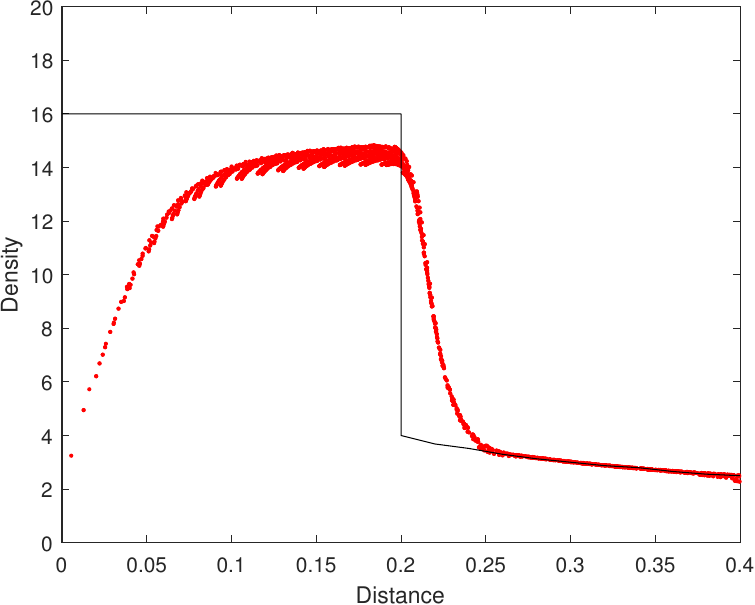}
    \caption{Numerical results of the Noh problem using the $Q^3-Q^2$ element pair with $h = \frac{1}{20}$ at $t = 0.6$: (a) density field, (b) density values compared with the analytical solution.}
    \label{Noh_Q3_Q2}
  \end{center}
\end{figure}

Figures~\ref{Noh_Q2_Q1} and~\ref{Noh_Q3_Q2} illustrate the numerical performance of the proposed method at $t = 0.6$. The right panels compare the numerical results with the analytical solution, with numbers indicating values at the degrees of freedom of the thermodynamic variable. For the same mesh, the number of DOFs in the $Q^3-Q^2$ element pair is approximately double that of the $Q^2-Q^1$ pair. Both element pairs perform well while maintaining satisfactory mesh quality. Without hourglass control, the Noh problem tends to produce slightly distorted meshes.

\subsection*{Dukowicz-Meltz} 
This problem is a piston-driven setup involving two regions filled with ideal gas. The left part of the domain is a right-angled trapezoid with a vertical left boundary, while the right part is a slanted parallelogram whose left boundary coincides with the right boundary of the trapezoid. The interface between the two regions is an oblique line inclined at $60^{\circ}$. A piston located at the left boundary compresses the gas in the left region with a velocity of $1.48$. The initial conditions for the left region are $\gamma = 1.4$, $\rho = 1$, and $e = 2.5$, while for the right region they are $\gamma = 1.4$, $\rho = 1.5$, and $e = 2.5$.

\begin{figure}
  \begin{center}
    \includegraphics[width=0.90\textwidth]{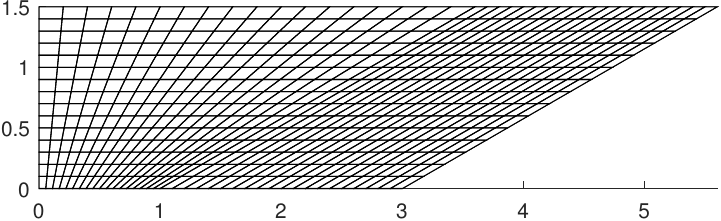}
    \caption{The initial mesh for Dukowicz-Meltz problem. }
    \label{DM_grid}
  \end{center}
\end{figure}

\begin{figure}
  \begin{center}
    \includegraphics[width=0.49\textwidth]{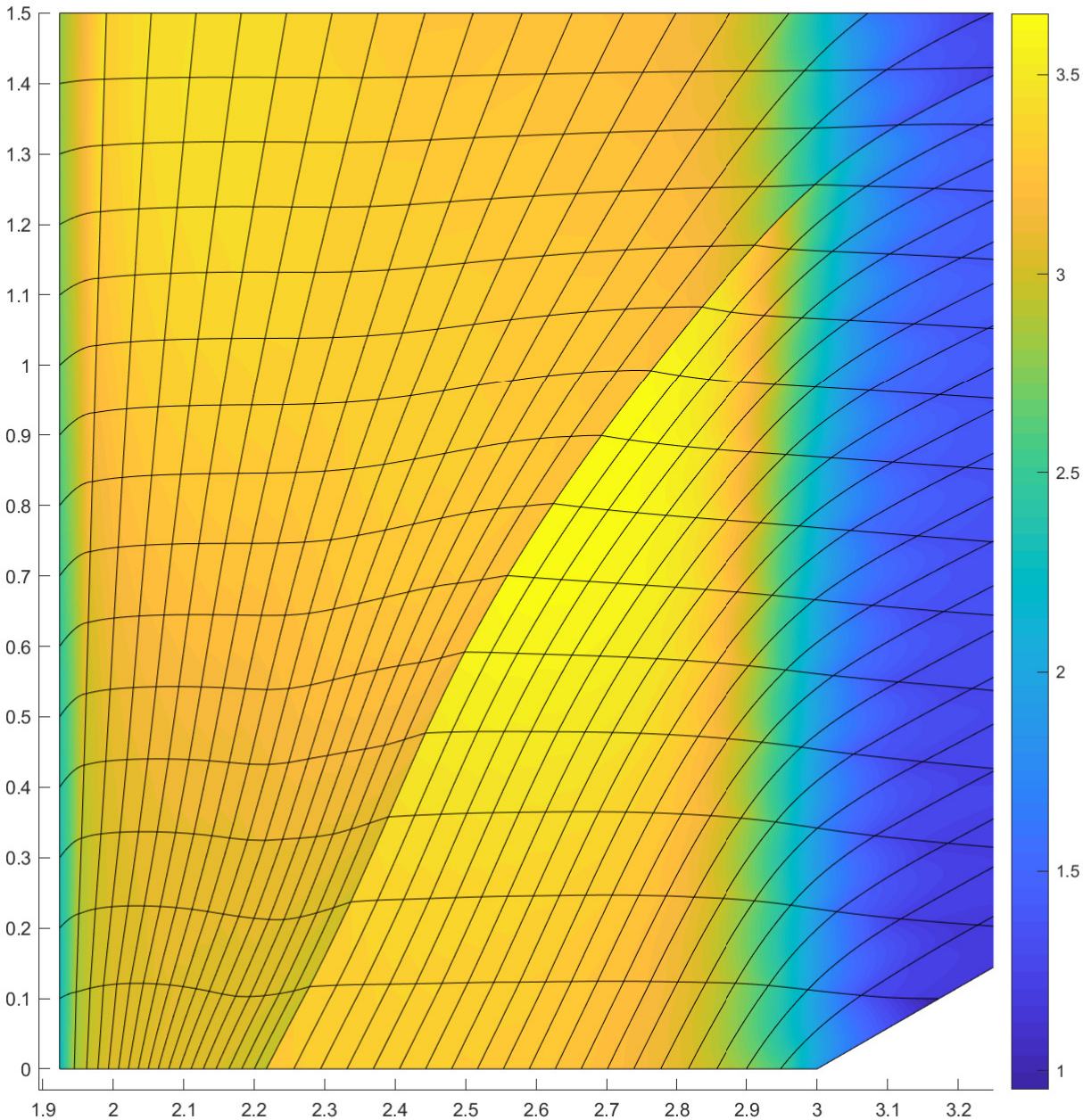}
    \includegraphics[width=0.49\textwidth]{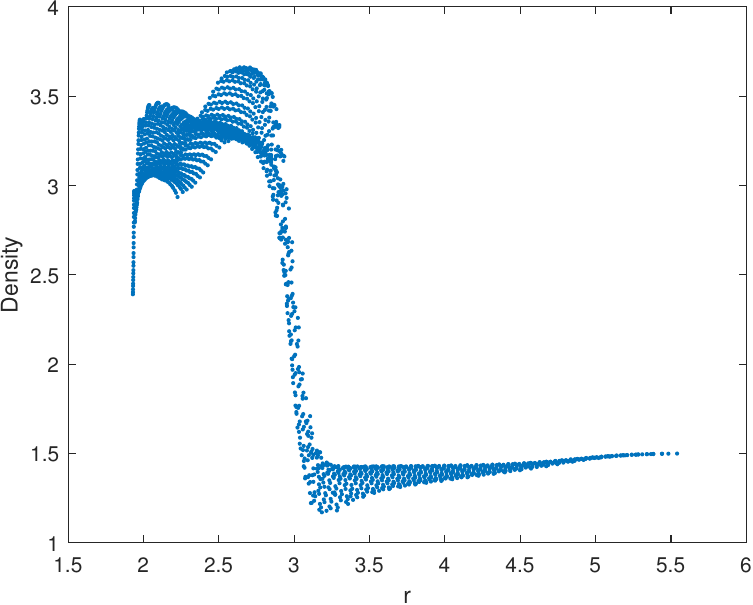}
    \caption{Numerical results of the Dukowicz-Meltz problem using the $Q^2-Q^1$ element pair at $t = 1.3$: (a) density field, (b) density values versus the $r$-coordinate.}
    \label{DM_Q2_Q1}
  \end{center}
\end{figure}

\begin{figure}
  \begin{center}
    \includegraphics[width=0.49\textwidth]{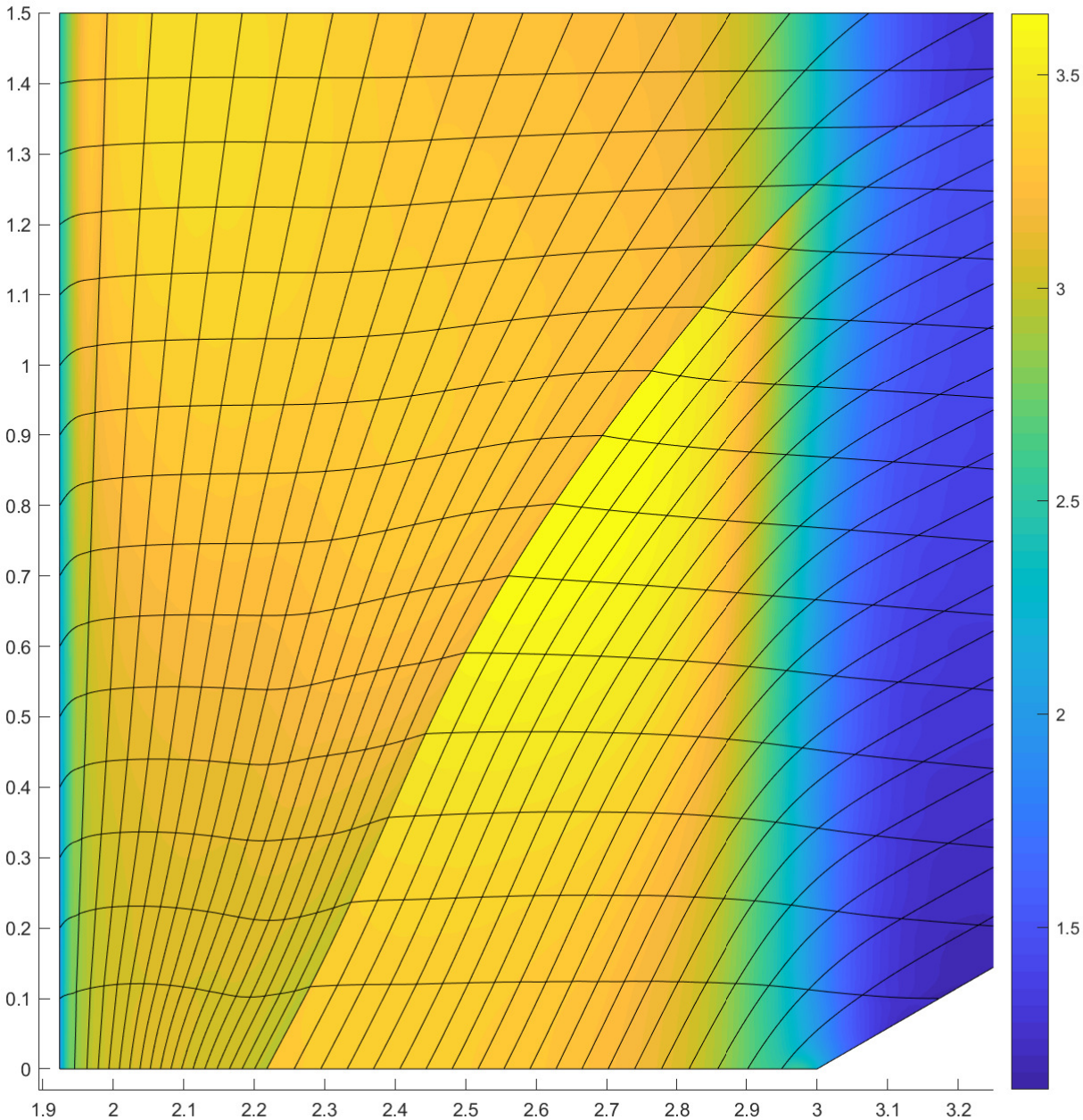}
    \includegraphics[width=0.49\textwidth]{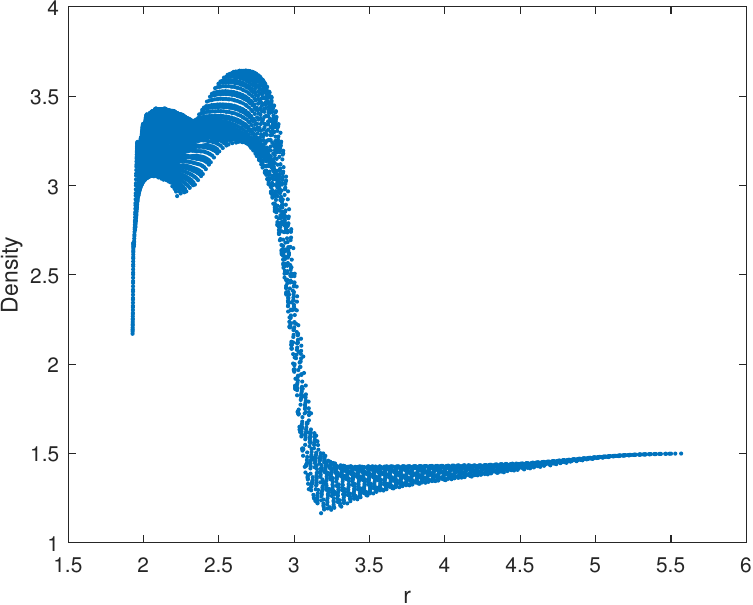}
    \caption{Numerical results of the Dukowicz-Meltz problem using the $Q^3-Q^2$ element pair at $t = 1.3$: (a) density field, (b) density values versus the $r$-coordinate.}
    \label{DM_Q3_Q2}
  \end{center}
\end{figure}

The initial mesh, shown in Figure~\ref{DM_grid}, consists of $38 \times 15$ elements. The left trapezoidal region is partitioned into $18 \times 15$ non-uniform elements, while the right parallelogram is uniformly partitioned into $20 \times 15$ elements. The terminal time is set to $t = 1.3$. Figures~\ref{DM_Q2_Q1} and~\ref{DM_Q3_Q2} present the numerical results for the $Q^2-Q^1$ and $Q^3-Q^2$ element pairs, respectively. The simulations demonstrate satisfactory numerical performance with good mesh quality.

\subsection*{Triple point}
The computational domain is $[0,7] \times [0,3]$ and consists of three distinct regions filled with ideal gases. The interactions among these gases generate a shock wave and a physical vortex. The left region, $[0,1] \times [0,3]$, has initial conditions $\gamma_1 = 1.5$, $\rho_1 = 1$, and $p_1 = 1$. The lower-right region, $[1,7] \times [0,1.5]$, is initialized with $\gamma_2 = 1.4$, $\rho_2 = 1$, and $p_2 = 0.1$. The upper-right region, $[1,7] \times [1.5,3]$, has $\gamma_3 = 1.6$, $\rho_3 = 0.125$, and $p_3 = 0.1$.

\begin{figure}
  \begin{center}
    \includegraphics[width=1.00\textwidth]{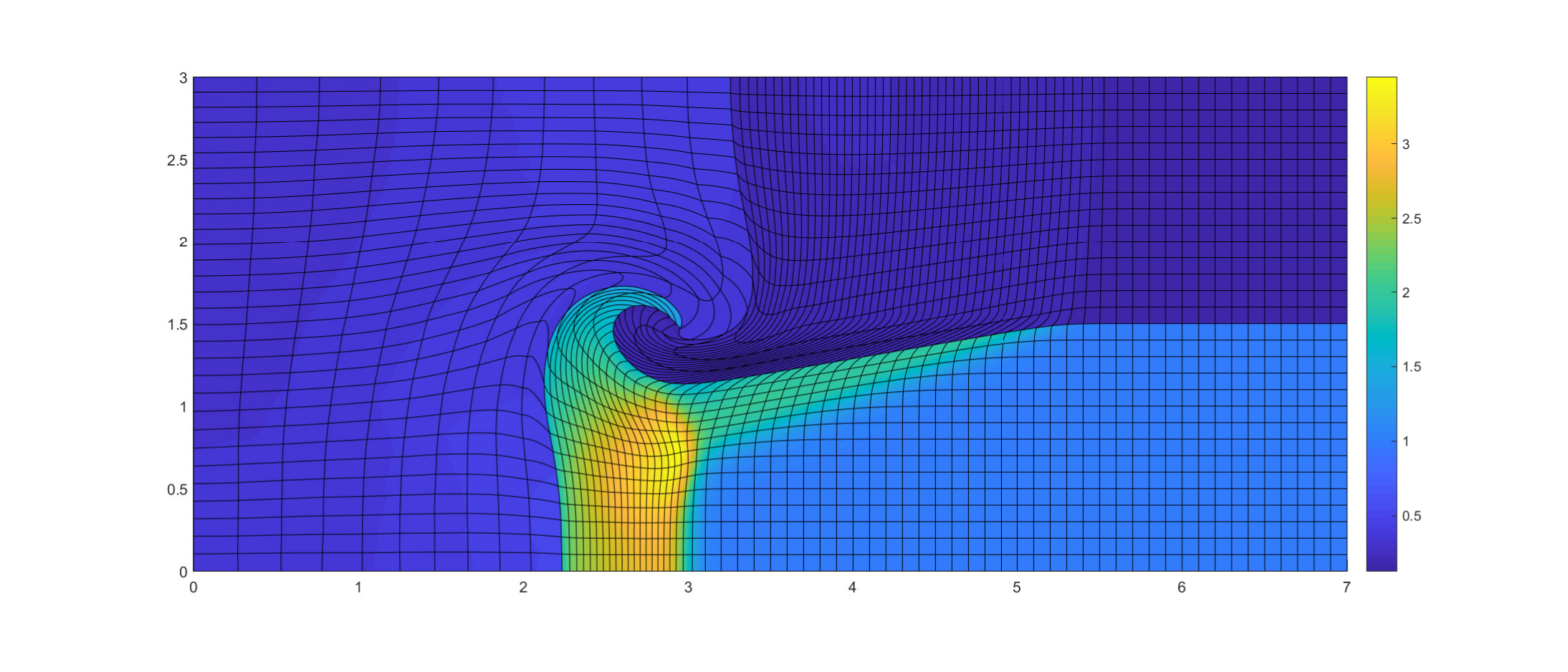}
    \caption{Density field for the triple-point problem using the $Q^2-Q^1$ element pair at $t = 2.5$.}
    \label{triple_1}
  \end{center}
\end{figure}

\begin{figure}
  \begin{center}
    \includegraphics[width=1.00\textwidth]{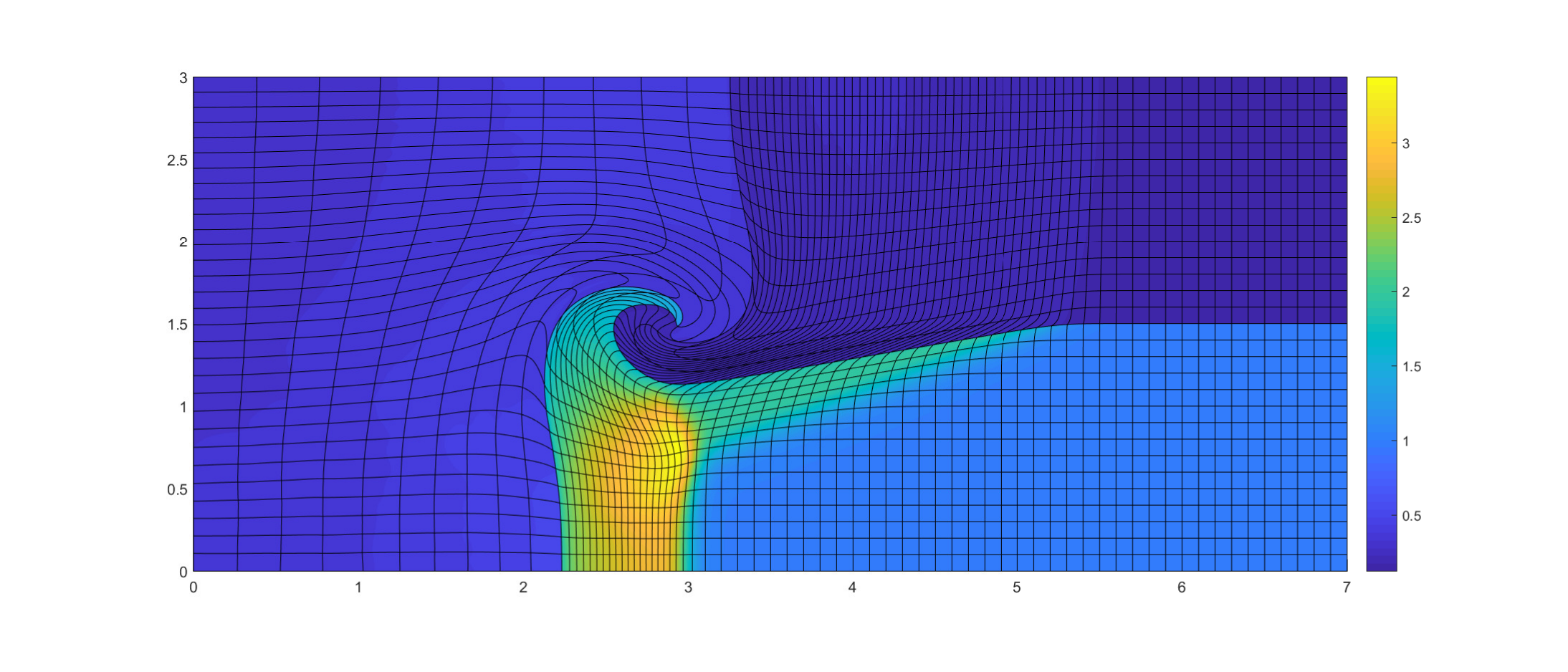}
    \caption{Density field for the triple-point problem using the $Q^3-Q^2$ element pair at $t = 2.5$.}
    \label{triple_2}
  \end{center}
\end{figure}

\begin{figure}
  \begin{center}
    \includegraphics[width=0.49\textwidth]{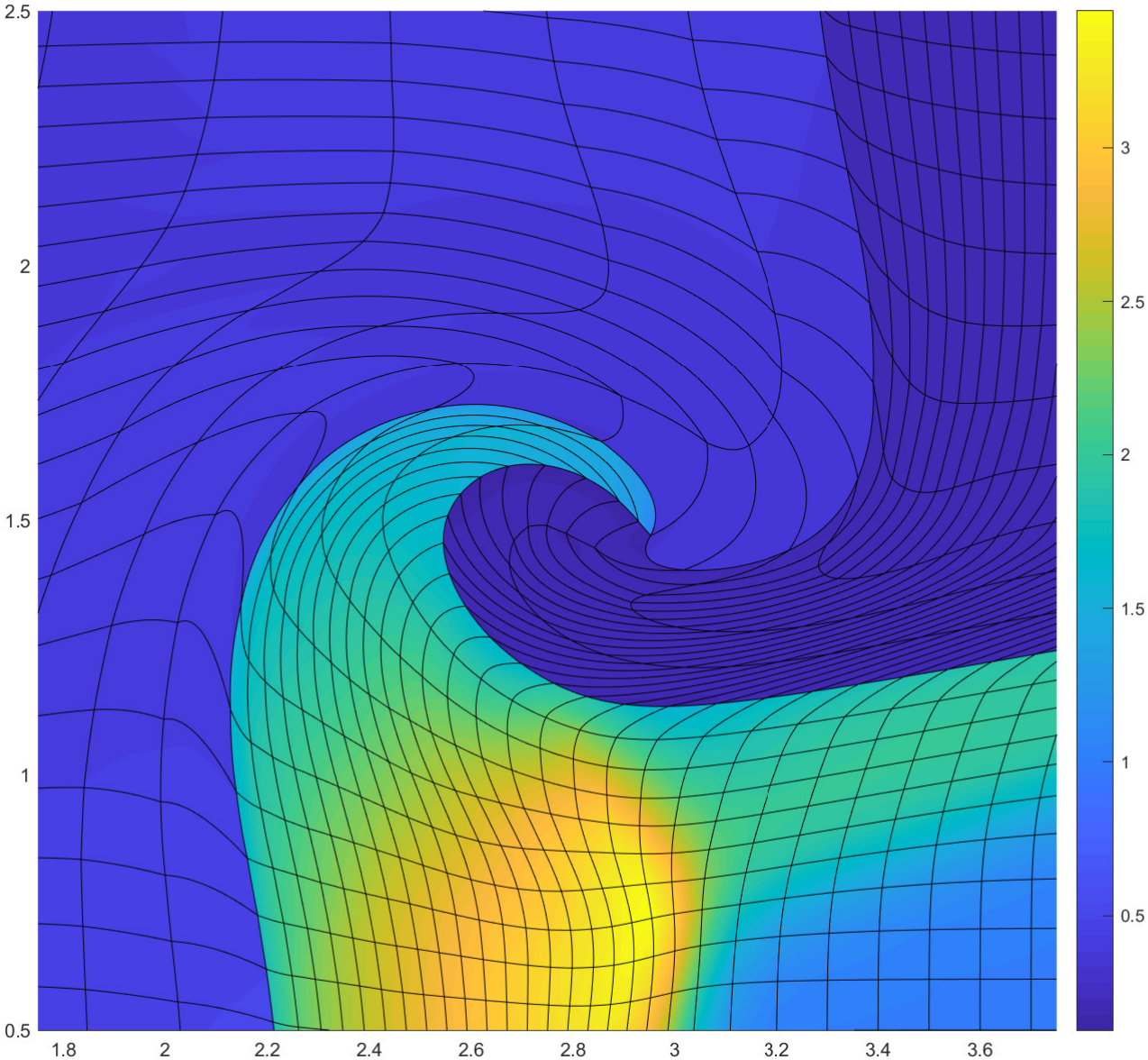}
    \includegraphics[width=0.49\textwidth]{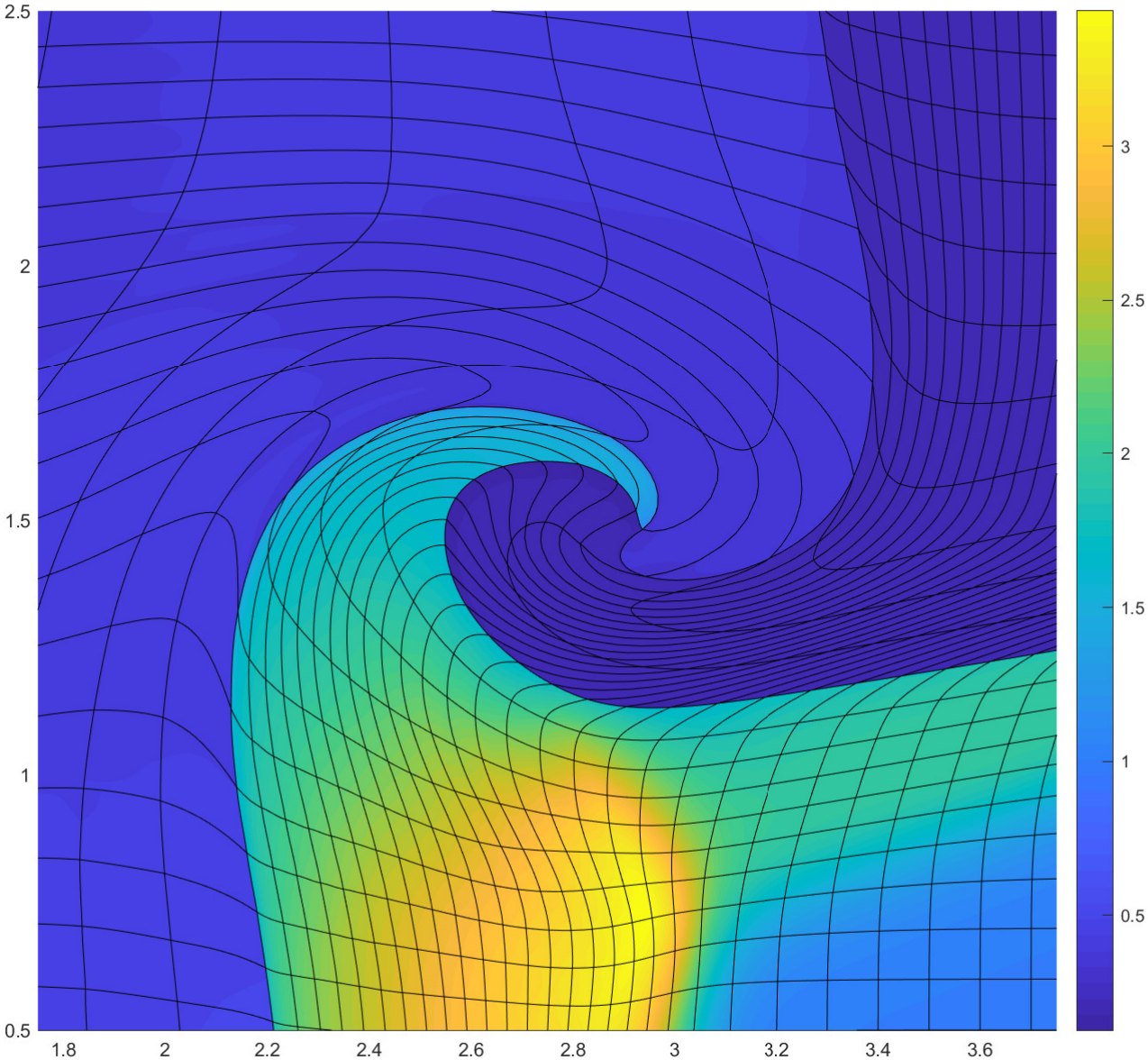}
    \caption{Zoomed-in density field for the triple-point problem at $t = 2.5$: (a) $Q^2-Q^1$ element pair, (b) $Q^3-Q^2$ element pair.}
    \label{triple_3}
  \end{center}
\end{figure}

Figures~\ref{triple_1} and~\ref{triple_2} show the numerical results for the entire computational domain using the $Q^2-Q^1$ and $Q^3-Q^2$ element pairs, respectively. Figure~\ref{triple_3} presents a zoomed-in view of the vortex region for both discretizations, highlighting that the $Q^3-Q^2$ solution captures more pronounced vortex twisting.

\subsection*{Sedov}
The purpose of presenting the numerical results of the Sedov problem is to provide an intuitive illustration of hourglass motion in the high-order case. The Sedov problem models an explosion under planar geometry, driven by a point source of internal energy deposited at the origin, where the resulting blast wave reaches a radius of 1 at $t = 1$. A stability issue affects the proposed high-order SGH: if the initial internal energy of the elements is set below $1 \times 10^{-3}$, or if the initial mesh is refined, the internal energy may become negative during the simulation. For the $Q^2-Q^1$ case, replacing the shape functions in equation~\eqref{eq:energy_hg_av_discrete} with the traditional $Q^1$ FEM shape functions restores stability. However, for the $Q^3-Q^2$ case, this approach is ineffective, and the initial internal energy is set to $1 \times 10^{-3}$.

\begin{figure}
  \begin{center}
    \includegraphics[width=0.49\textwidth]{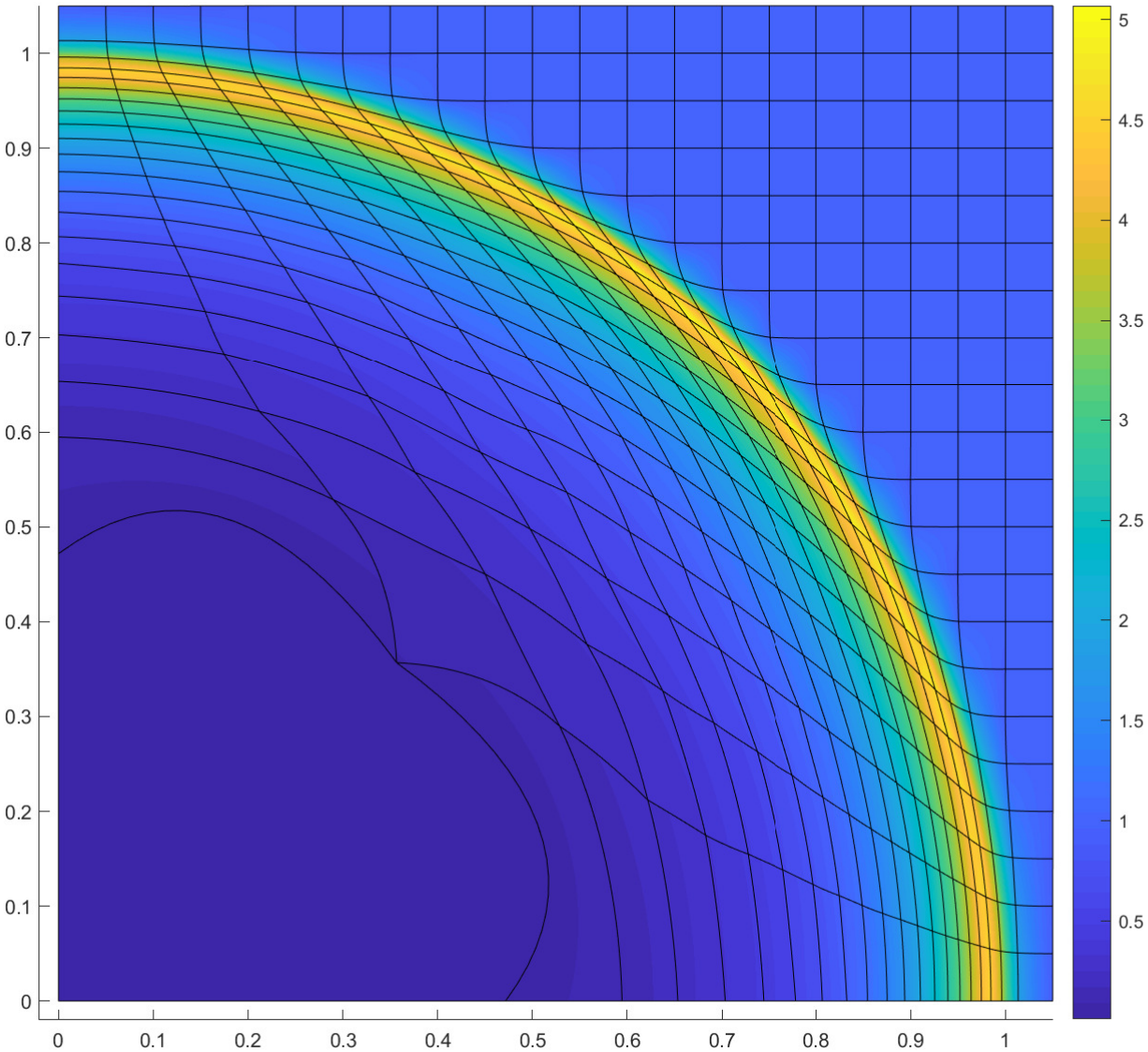}
    \includegraphics[width=0.49\textwidth]{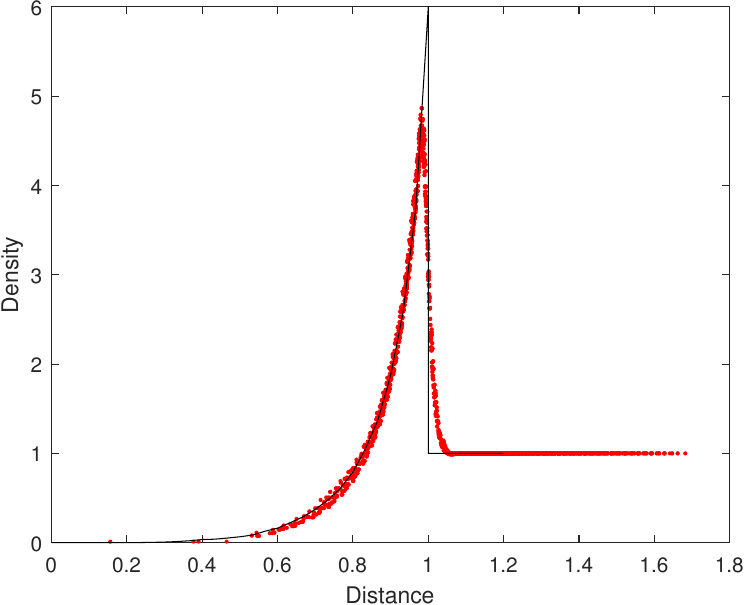}
    \caption{Numerical results for the Sedov problem using the $Q^2-Q^1$ element pair at $t = 1$ without hourglass control. (a) mesh distribution and density field .(b) density values compared with the analytical solution.}
    \label{Sedov_nohg}
  \end{center}
\end{figure}

\begin{figure}
  \begin{center}
    \includegraphics[width=0.49\textwidth]{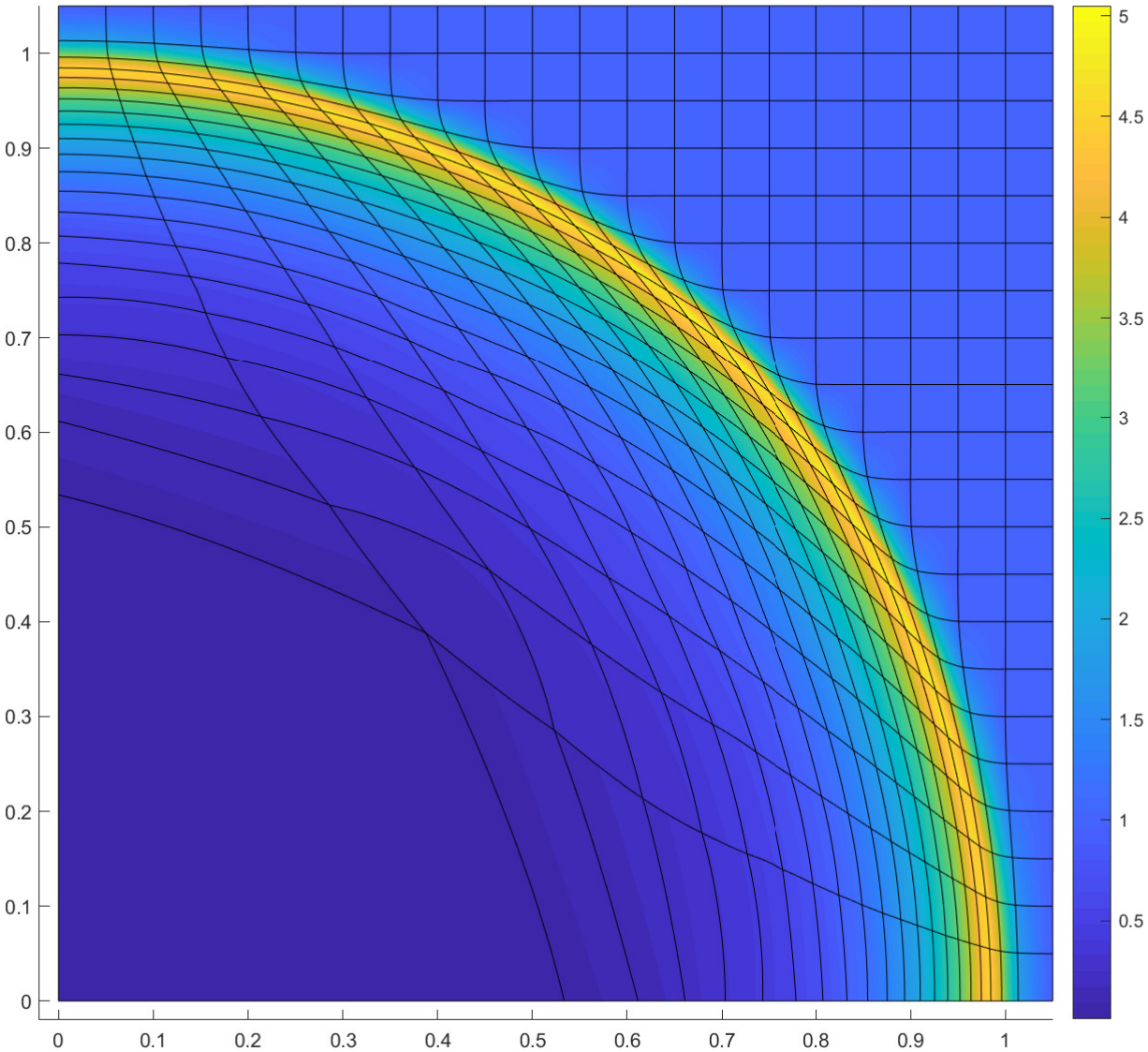}
    \includegraphics[width=0.49\textwidth]{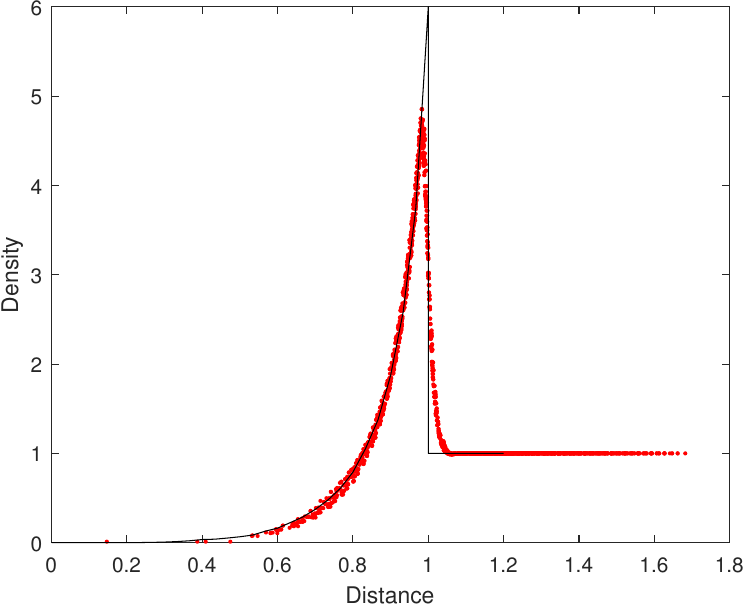}
    \caption{Numerical results for the Sedov problem using the $Q^2-Q^1$ element pair at $t = 1$ with hourglass control. (a) mesh distribution and density field .(b) density values compared with the analytical solution.}
    \label{Sedov_hg}
  \end{center}
\end{figure}

\begin{figure}
  \begin{center}
    \includegraphics[width=0.49\textwidth]{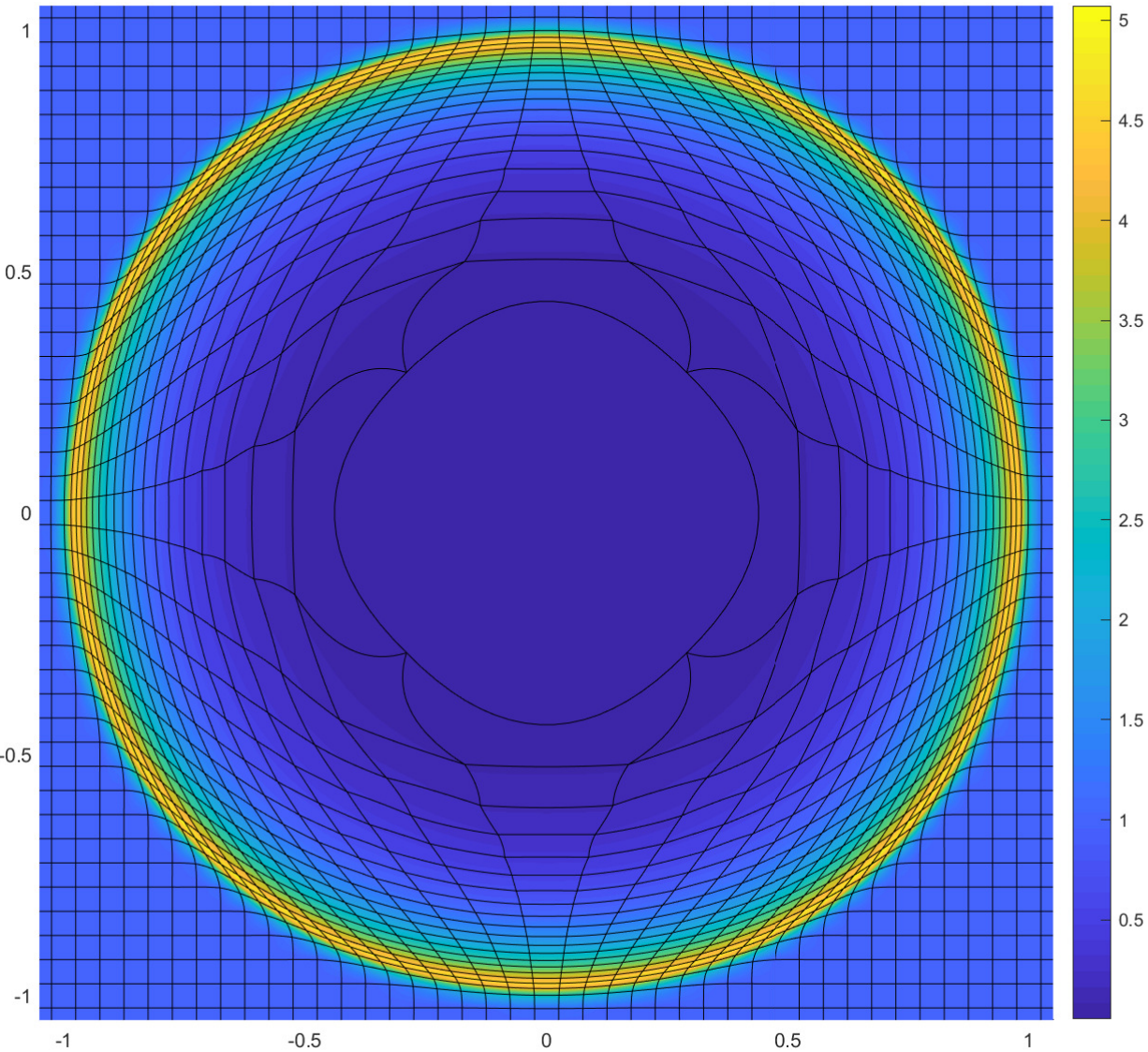}
    \includegraphics[width=0.49\textwidth]{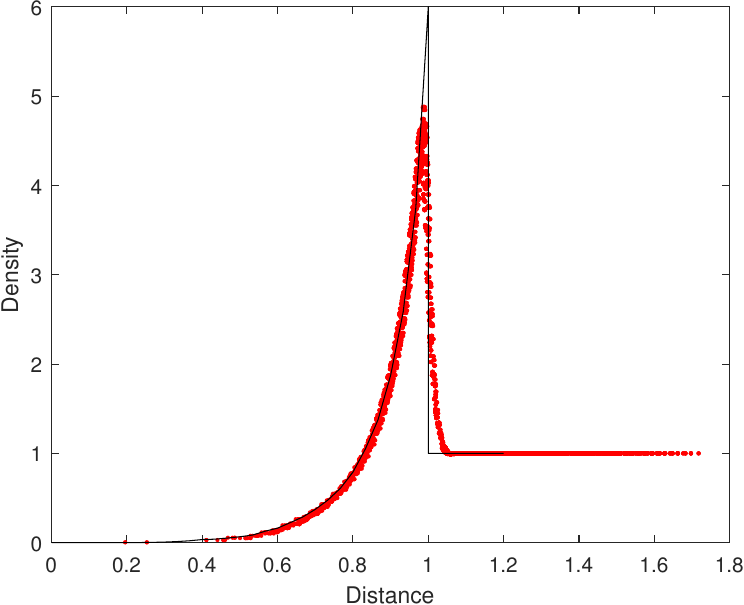}
    \caption{Numerical results for the Sedov problem(center point source) using the $Q^2-Q^1$ element pair at $t = 1$ without hourglass control. (a) mesh distribution and density field .(b) density values compared with the analytical solution.}
    \label{Source_nohg}
  \end{center}
\end{figure}

\begin{figure}
  \begin{center}
    \includegraphics[width=0.49\textwidth]{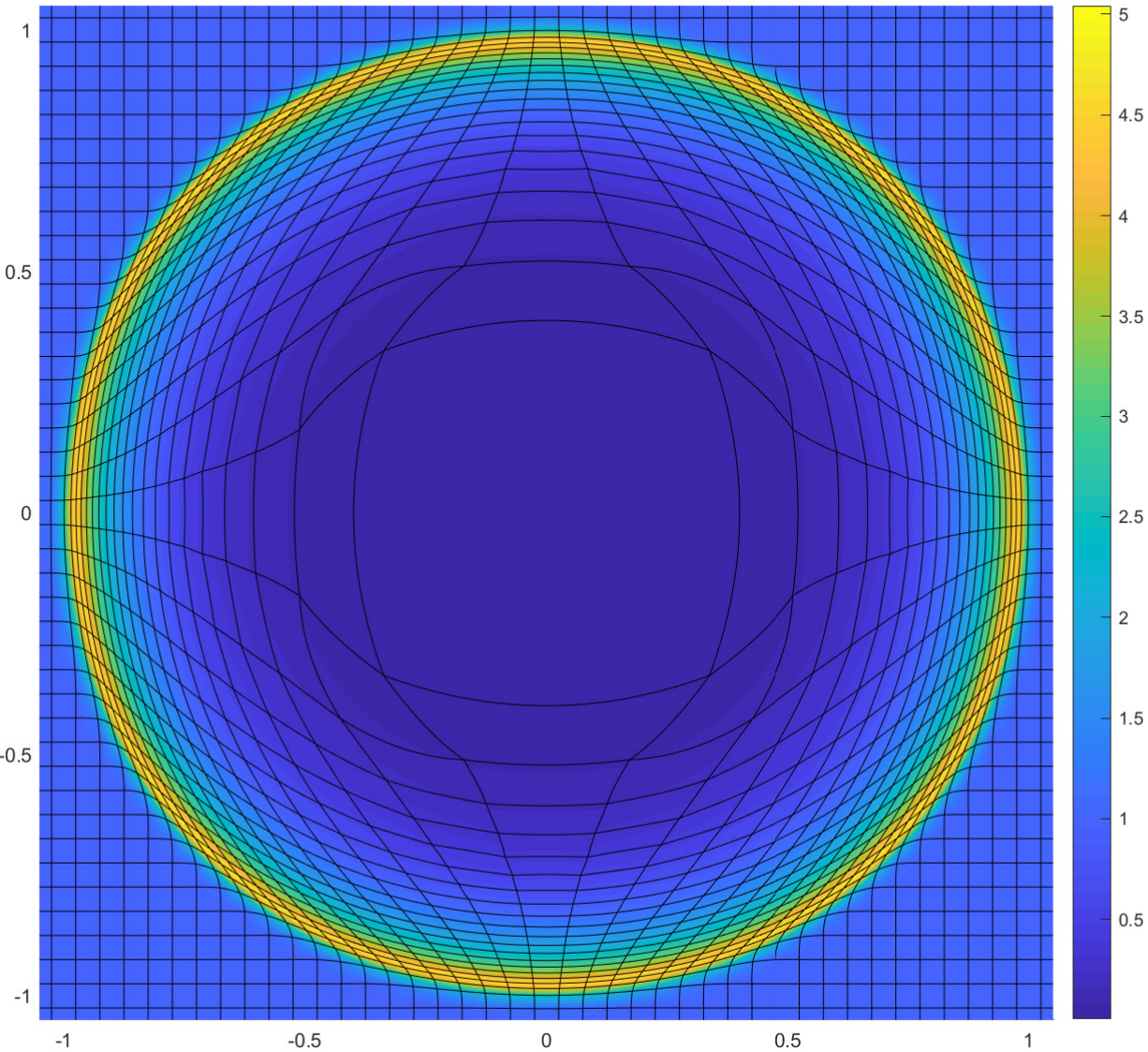}
    \includegraphics[width=0.49\textwidth]{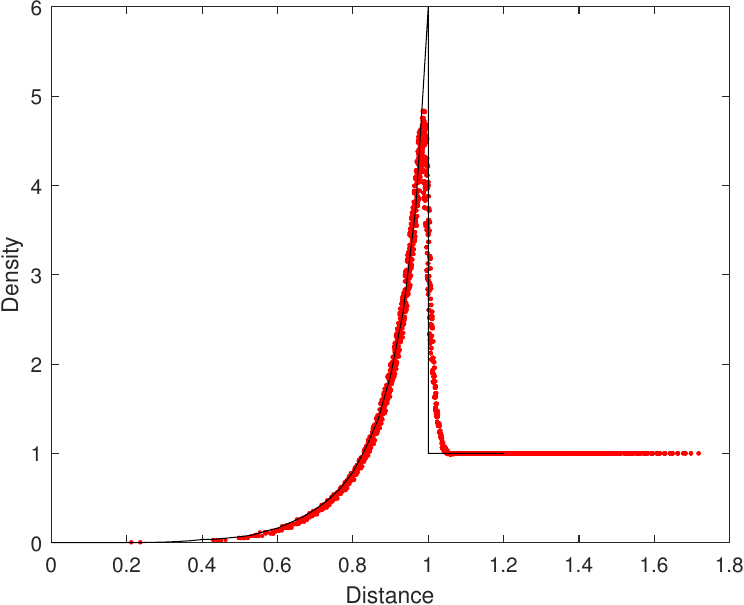}
    \caption{Numerical results for the Sedov problem(center point source) using the $Q^2-Q^1$ element pair at $t = 1$ with hourglass control. (a) mesh distribution and density field .(b) density values compared with the analytical solution.}
    \label{Source_hg}
  \end{center}
\end{figure}

Figure~\ref{Sedov_nohg} shows the $Q^2-Q^1$ results with a point source at the origin without hourglass control, while Figure~\ref{Sedov_hg} presents the corresponding results with hourglass control. Figure~\ref{Source_nohg} depicts the no hourglass control results with a point source located at the central elements, and Figure~\ref{Source_hg} shows the results with hourglass control.

Because artificial viscosity suppresses hourglass motion, the grid quality near the shock remains good. However, away from the shock, where the viscosity is inactive, hourglass distortions appear in Figures~\ref{Sedov_nohg} and~\ref{Source_nohg}. Despite this, the density field values remain nearly unchanged. Although the poor mesh quality does not immediately affect accuracy, it may lead to instability and chaotic behavior in subsequent simulations.

\begin{figure}
  \begin{center}
    \includegraphics[width=0.49\textwidth]{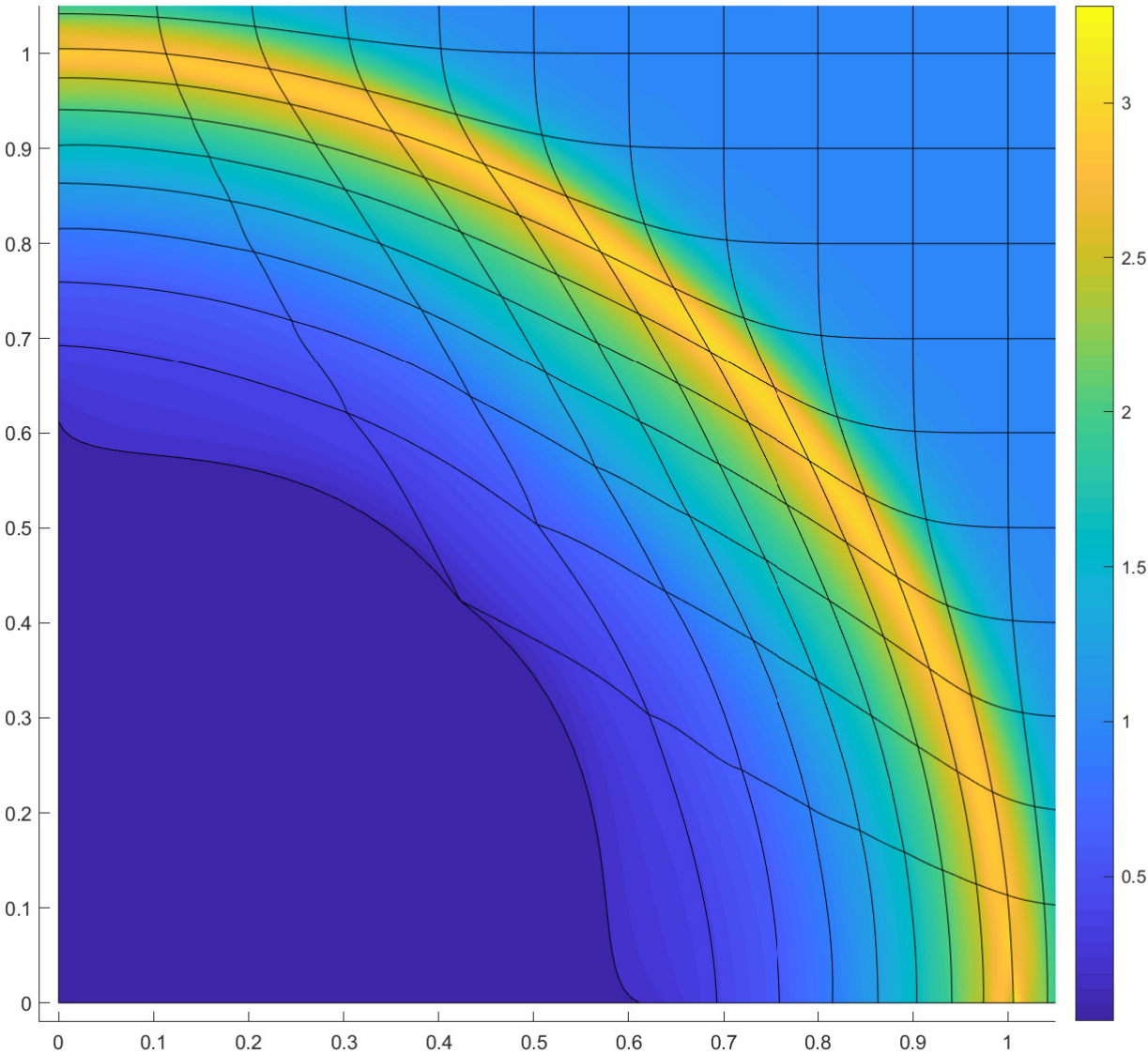}
    \includegraphics[width=0.49\textwidth]{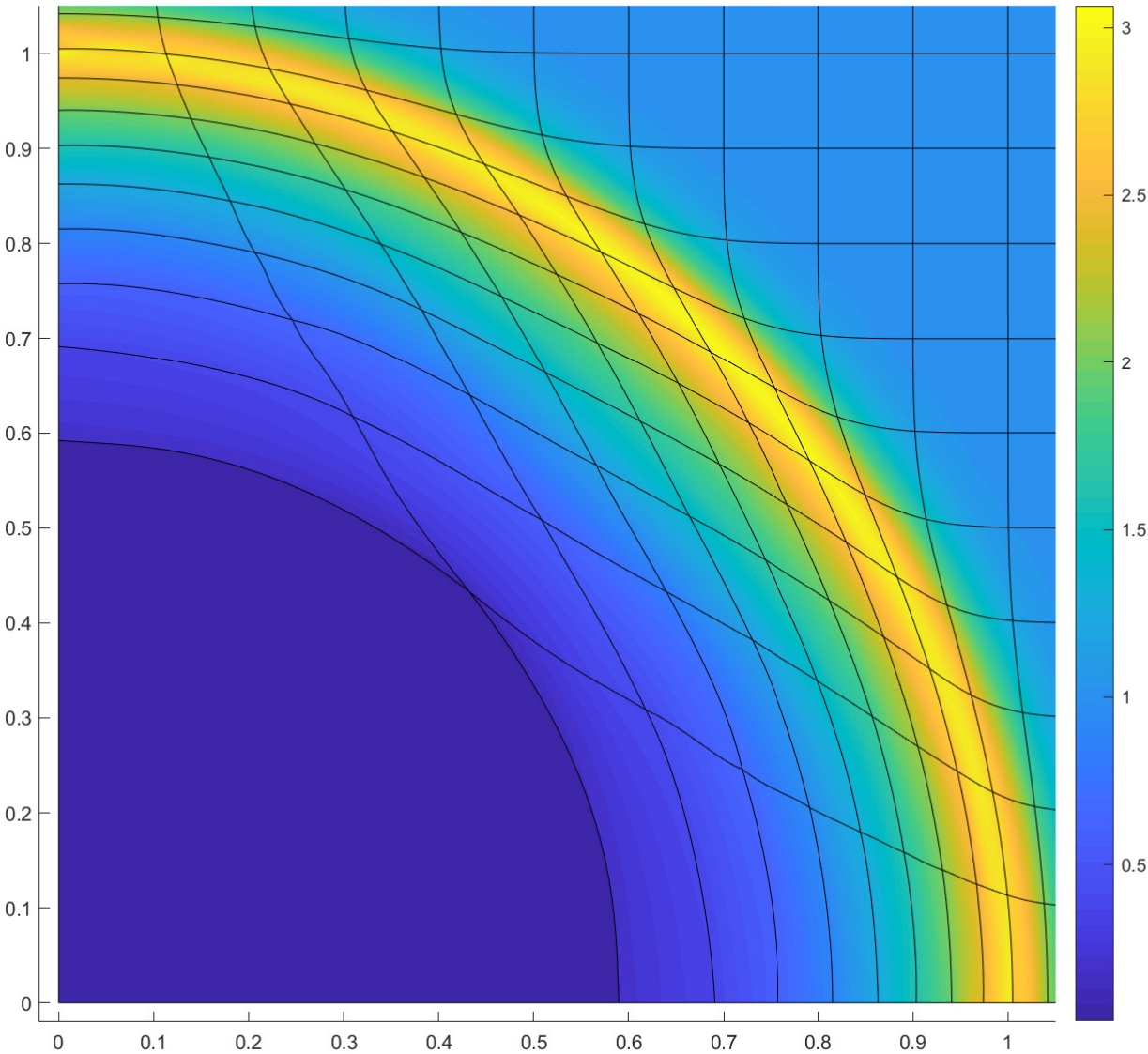}
    \caption{Density field and mesh distribution for the Sedov problem using the $Q^3-Q^2$ element pair at $t = 1$.(a) without hourglass control.(b) with hourglass control.}
    \label{Sedov_Q3}
  \end{center}
\end{figure}

\begin{figure}
  \begin{center}
    \includegraphics[width=0.49\textwidth]{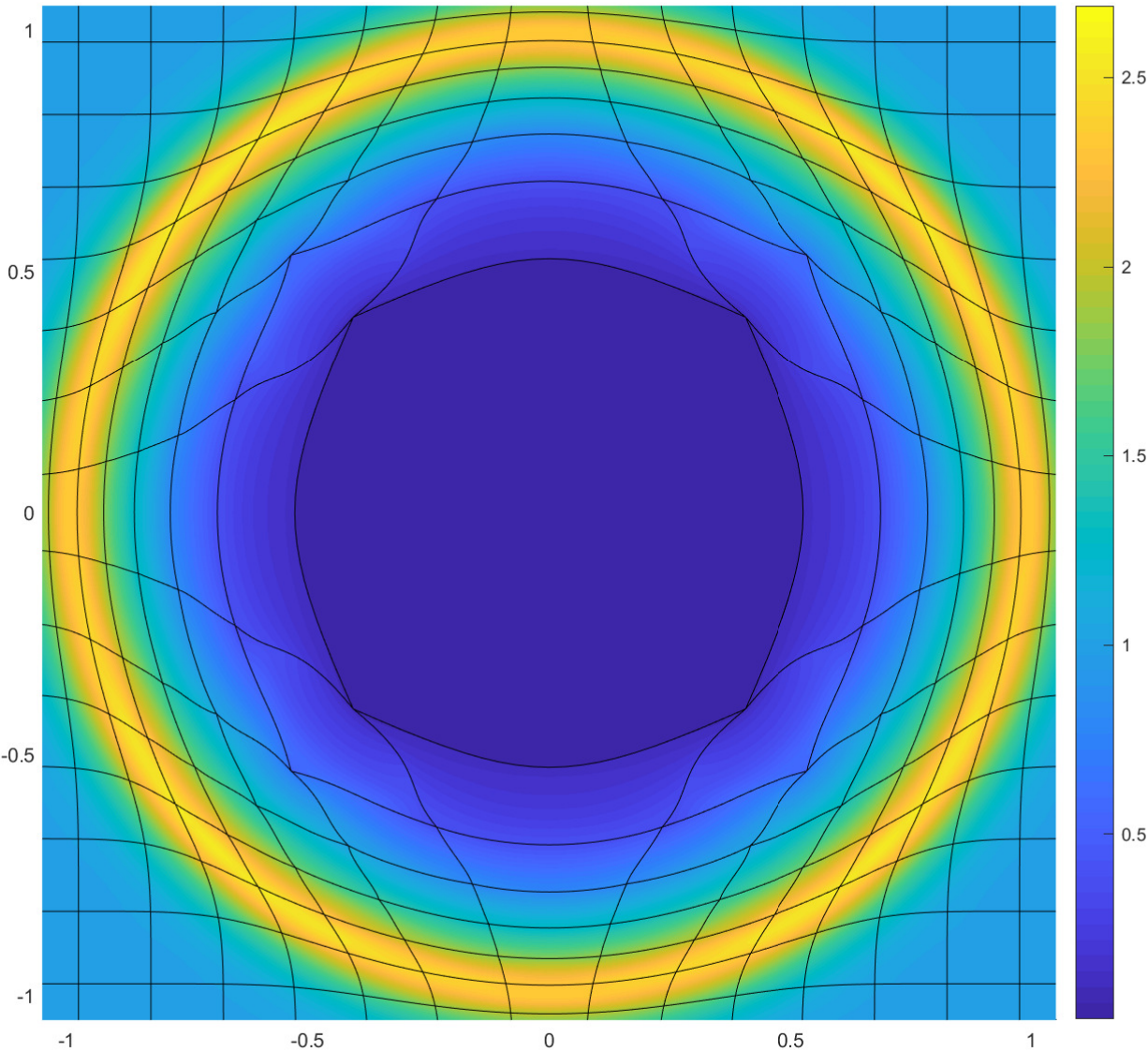}
    \includegraphics[width=0.49\textwidth]{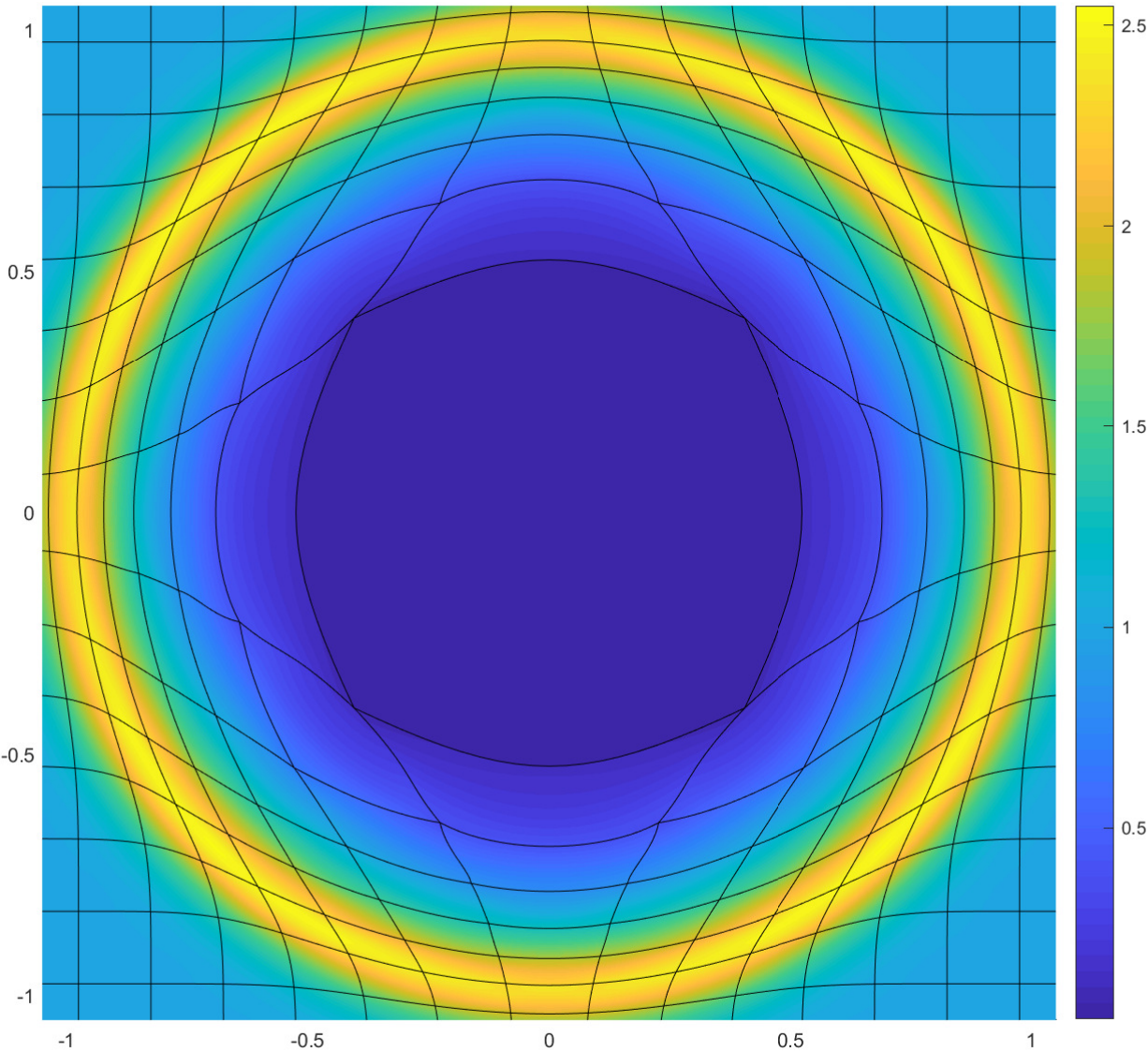}
    \caption{Density field and mesh distribution for the Sedov problem(center point source) using the $Q^3-Q^2$ element pair at $t = 1$.(a) without hourglass control.(b) with hourglass control.}
    \label{Source_Q3}
  \end{center}
\end{figure}

Figures~\ref{Sedov_Q3} and~\ref{Source_Q3} illustrate the same scenario for the $Q^3-Q^2$ case. Hourglass distortions appear in the mesh regions away from the shock. Due to stability constraints, the initial internal energy is set to $1 \times 10^{-3}$, which is not physically realistic; therefore, the density values versus the radius are not shown.

\section{Conclusion}\label{sec:conclusion}
This work completes the framework of high-order staggered Lagrangian hydrodynamics (SGH) with artificial viscosity and hourglass control. For the hourglass control algorithm, we extended the subzonal pressure hourglass control method to high-order case. From an interpolation perspective, the anti-hourglass force was constructed based on pressure variations, resulting in a unified framework for this class of hourglass control algorithms. Although no novel artificial viscosity was introduced, we highlighted the special structure and provided explicit expressions for the prerequisite variables required in its computation. Our efforts focused on reducing computational cost and enabling efficient implementation. Both algorithms require numerical integration of one order higher than that used for the conservation law discretization. This choice, however, reduces the robustness of the proposed SGH scheme, indicating that robustness and efficiency cannot be fully achieved simultaneously. For smooth, shock-free problems, the observed convergence rates match theoretical predictions, and high-order schemes demonstrate superior efficiency per degree of freedom, achieving the same accuracy with fewer DOFs. For problems involving shocks, the proposed high-order SGH scheme produces high-fidelity results when appropriate artificial viscosity and hourglass control parameters are employed.


\section*{Acknowledgment}
The research is supported by the National Natural Science Foundation 
of China(Grant No. 12571414). 


\bibliographystyle{abbrv}
\bibliography{ref}
\end{document}